\documentclass[10pt,a4paper,oneside]{article}
\usepackage[utf8]{inputenc}
\usepackage[english]{babel}
\pdfoutput=1

\usepackage{amsmath,amsfonts,amssymb,amsopn,amscd,amsthm}
\usepackage{mathtools}
\allowdisplaybreaks

\usepackage[yyyymmdd,hhmmss]{datetime}

\usepackage{fancyhdr}
\usepackage{amsthm}

\newtheorem{theorem}{Theorem}[section]

\newtheorem{proposition}[theorem]{Proposition}

\newtheorem{definition}[theorem]{Definition}

\usepackage{amsmath,amsfonts,amssymb,amsopn,amscd}
\usepackage{mathtools}

\newcommand{\RR}{\mathbb{R}}
\newcommand{\RRplus}{\mathbb{R}_+}
\newcommand{\NN}{\mathbb{N}}
\newcommand{\ZZ}{\mathbb{Z}}

\newcommand{\DescSet}[2]{\{#1\,|\,#2\}}
\newcommand{\EnumSet}[1]{\{#1\}}
\newcommand{\LRDescSet}[2]{\left\{#1\,\middle|\,#2\right\}}

\newcommand{\Iverson}[1]{[#1]}

\newcommand{\Modulus}[1]{|#1|}
\newcommand{\LRModulus}[1]{\left|#1\right|}
\newcommand{\BigModulus}[1]{\Bigl|#1\Bigr|}
\newcommand{\Cardinality}[1]{|#1|}
\newcommand{\SymDiff}{\bigtriangleup}

\newcommand{\D}[1]{{\mathrm{d}\,\!#1}}
\newcommand{\Expect}{{\mathbb{E}}}

\newcommand{\Dimension}{d}
\newcommand{\RRD}{\RR^{\Dimension}}
\newcommand{\RRplusD}{\RRplus^{\Dimension}}
\newcommand{\VolumeCoeff}[1]{v_{#1}}
\newcommand{\Distance}[1]{||#1||}
\newcommand{\Diameter}[1]{\operatorname{diam}(#1)}
\newcommand{\SetDistance}[1]{{\delta(#1)}}
\newcommand{\Lebesgue}{\mathcal{L}}
\newcommand{\Radius}{R}
\newcommand{\Borel}{\mathcal{B}}
\newcommand{\BoundedBorel}{\mathcal{B}_b}
\newcommand{\Sphere}[1]{\mathcal{S}(#1)}
\newcommand{\Ring}[1]{\mathcal{R}(#1)}
\newcommand{\Complement}[1]{{#1^{c}}}
\newcommand{\Configs}[1]{\mathcal{C}_{#1}}
\newcommand{\Algebra}[1]{\mathcal{F}_{#1}}
\newcommand{\ProductAlgebra}[2]{\mathcal{F}_{#1}^{\otimes{#2}}}

\newcommand{\RConnectedToIn}[1]{
 {\,\xleftrightarrow[]{\text{in }{#1}}\,}
}

\newcommand{\Law}{\mathcal{P}}
\newcommand{\Poisson}[1]{\Law^{\text{poi}}_{#1}}
\newcommand{\CritValuePercolation}[1]{\lambda_{b}(#1)}

\newcommand{\SymbolIsHS}{\mathcal{H}}
\newcommand{\IsHS}[1]{\SymbolIsHS(#1)}
\newcommand{\Potential}{u}
\newcommand{\ReducedPairCorrelation}{\rho}
\newcommand{\SymbolPartFun}{\mathcal{Z}}
\newcommand{\PartFun}[1]{\SymbolPartFun(#1)}
\newcommand{\HardSphere}[1]{\Law^{\text{hs}}_{#1}}
\newcommand{\GibbsMeasures}[1]{\mathcal{G}_{#1}}
\newcommand{\CERadius}[1]{\lambda_{ce}(#1)}

\newcommand{\Order}{\prec}
\newcommand{\LT}{\prec}
\newcommand{\LTE}{\preceq}

\newcommand{\HSPoiThin}[1]{\Law^{\text{thin}}_{#1}}
\newcommand{\OneShot}[1]{{p_{#1}}}
\newcommand{\ThinProba}[2]{{p(#1\mid{}#2)}}
\newcommand{\ThinChoice}[1]{{c(#1)}}

\newcommand{\TwistedZone}[1]{\Law^{\text{tw-zone}}_{#1}}
\newcommand{\TwistedRec}[1]{\Law^{\text{tw-rec}}_{#1}}
\newcommand{\DisagreementConnection}[2]{{\mathcal{D}_{#1}^{#2}}}

\newcommand{\SymbolHSSize}{\mathfrak{s}}
\newcommand{\HSSize}[1]{\SymbolHSSize(#1)}

\allowdisplaybreaks

\newcommand{\Support}{
The author acknowledges the support of the VIDI project ``Phase transitions, Euclidean fields and random fractals'', NWO 639.032.916.
}

\newcommand{\TitleFull}{Disagreement percolation for the hard-sphere model}

\newcommand{\TitleShort}{Hard-sphere\\ disagreement percolation}

\newcommand{\AuthorsFull}{Hofer-Temmel Christoph\footnote{\Email}\thanks{\Support}}

\newcommand{\AuthorsShort}{Hofer-Temmel}

\newcommand{\Email}{math@temmel.me}

\newcommand{\Keywords}{Keywords:
hard-sphere model,
disagreement percolation,
unique Gibbs measure,
stochastic domination,
Boolean model,
absence of phase transition,
dependent thinning
}

\newcommand{\Head}{
 \maketitle
 
\begin{abstract}
Disagreement percolation connects a Gibbs lattice gas and i.i.d. site percolation on the same lattice such that non-percolation implies uniqueness of the Gibbs measure.
This work generalises disagreement percolation to the hard-sphere model and the Boolean model.
Non-percolation of the Boolean model implies the uniqueness of the Gibbs measure and exponential decay of pair correlations and finite volume errors.
Hence, lower bounds on the critical intensity for percolation of the Boolean model imply lower bounds on the critical activity for a (potential) phase transition.
These lower bounds improve upon known bounds obtained by cluster expansion techniques.
The proof uses a novel dependent thinning from a Poisson point process to the hard-sphere model, with the thinning probability related to a derivative of the free energy.

\end{abstract}
{}
 \Keywords{}\\\MSC{}
}

\title{\TitleFull{}}
\author{\AuthorsFull{}}
\date{}
\begin{document}

\Head{}
\section{Introduction}

\par
Disagreement percolation by van den Berg and Maes~\cite{VanDenBerg_Maes__DisagreementPercolationInTheStudyOfMarkovFields__AP_1994} is a sufficient condition on the activity of a discrete Gibbs specification on a graph for uniqueness of the Gibbs measure.
It implies the absence of phase transitions and the analyticity of the free energy in the high-temperature case.
It has also been used to derive the Poincar\'{e} inequality in the context of lattice Ising spin systems~\cite{Chazottes_Redig_Voellering__ThePoincareInequalityForMarkovRandomFieldsProvedViaDisagreementPercolation__IMNS_KNAW_2011}.
This paper generalises disagreement percolation to the hard-sphere model on $\RRD$, the continuum equivalent of the well-studied hard-core model~\cite{VanDenBerg_Steif__PercolationAndTheHardCoreLatticeGasModel__SPA_1994}.

\par
The core of disagreement percolation is a coupling between three point processes on a bounded domain.
Two are hard-sphere models with the same activity and differing boundary conditions.
The third one is a Boolean model stochastically dominating the points of disagreement between the two hard-sphere models.
The connected components of the Gilbert graph of the Boolean model connected to the boundary control the extent of the differing influence of the boundary conditions on the hard-sphere models.
In the sub-critical phase of percolation, the almost-sure finite percolation clusters imply the equality of the two hard-sphere realisations with high probability on a small domain inside a larger domain.
Taking a limit along an exhaustive sequence of bounded domains implies the uniqueness of the Gibbs measure of the hard-sphere model.

\par
The disagreement coupling connects the activity of the hard-sphere models and the intensity of the Poisson point process.
Hence, lower bounds on the critical intensity of the Boolean model imply lower bounds on the critical activity of the hard-sphere model.
In one dimension, the results replicate Tonk's classic result of the complete absence of phase transitions~\cite{Tonks__TheCompleteEquationOfStateOfOneTwoAndThreeDimensionalGasesOfElasticHardSpheres__PhysRev_1936}.
In two dimensions, the new bounds improve upon the best known cluster expansion bounds~\cite{Ruelle__StatisticalMechanics_RigorousResults__Ben_1969,Fernandez_Procacci_Scoppola__TheAnalyticityRegionOfTheHardSphereGas_ImprovedBounds__JSP_2007} by at least a factor of two.
Scaling discrete results suggest that they exceed the best theoretical largest activities achievable by cluster expansion techniques.
In high dimensions, extrapolation of known upper bounds on the activities achievable in the discrete case to the continuum suggests that the disagreement percolation bounds always go beyond the region attainable by cluster expansion techniques.

\par
This work exclusively treats the hard-sphere model.
One reason is its central importance in statistical mechanics and its easy and emblematic definition.
Another reason is the comparison with the cluster expansion bounds.
More important though, the bounds in this paper stem from a twisted disagreement coupling optimised for the hard-sphere model.
While a generalisation of the disagreement approach to simple finite-range Gibbs point processes with bounded interaction range seems possible, the twisted coupling depends critically on the hard-sphere constraint.
The twisted approach takes inspiration from a disagreement percolation tailored to the hard-core model~\cite{VanDenBerg_Steif__PercolationAndTheHardCoreLatticeGasModel__SPA_1994}.

\par
The twisted disagreement coupling is defined in a recursive fashion and uses conditional couplings between a hard-sphere model and its dominating Poisson point process.
The measurability of such a conditional coupling with respect to its boundary conditions is crucial for the existence of the twisted disagreement coupling.
The measurability of dominating couplings has not been a topic in the relevant literature on couplings~\cite{Preston__SpatialBirthAndDeathProcesses__BIIS_1975,Georgii_Kueneth__StochasticComparisonOfPointRandomFields__JAP_1997} yet.
One solution is the use a dependent thinning from the dominating Poisson point process.
The thinning probability is the derivative of the free energy of the yet unexplored part of the domain, rescaled by the activity.
It can be expressed as a ratio of partition functions.
The thinning approach is the key to ignore the uncountable nature of $\RRD$ and to focus on the almost-surely finite set of points of interest.

\par
Section~\ref{sec_setup} introduces notation and basic terms.
The main theorems, resulting bounds and discussion are in Section~\ref{sec_results}.
Section~\ref{sec_daperc_proof} contains the proofs about disagreement percolation.
Section~\ref{sec_thin} presents the dependent thinning.
Section~\ref{sec_twisted} elaborates the twisted disagreement coupling.

\section{Setup}
\label{sec_setup}

\subsection{Space}
\label{sec_space}

\par
Consider the \emph{Euclidean space} $\RRD$ with the \emph{Euclidean metric} $\Distance{.}$ and the \emph{Lebesgue measure} $\Lebesgue$.
The \emph{bounded} and \emph{all Borel sets} of $\RRD$ are $\BoundedBorel$ and $\Borel$ respectively.
Fix a \emph{non-negative finite radius} $\Radius$.
For $x\in\RRD$, let $\Sphere{x}:=\DescSet{y\in\RRD}{\Distance{x-y}\le\Radius}$ be the \emph{closed sphere} of radius $\Radius$ around $x$.
The \emph{volume} of $\Sphere{x}$ is $\VolumeCoeff{\Dimension}\Radius^\Dimension$.
For $B\in\Borel$, let $\Sphere{B}:=\bigcup_{x\in{}B}\Sphere{x}$ and $\Ring{B}:=\Sphere{B}\setminus{}B$ be the \emph{sphere} and \emph{ring} of radius $\Radius$ around $B$ respectively.
Let $\SetDistance{A,B}$ be the \emph{distance} between $A,B\in\BoundedBorel$.
A \emph{van~Hove sequence}~\cite[Def~2.1.1]{Ruelle__StatisticalMechanics_RigorousResults__Ben_1969} is a monotone increasing sequence $(B_n)_{n\in\NN}$ of bounded Borel sets converging to $\RRD$ and eventually containing every bounded Borel set.
The increasing hypercubes $([-n,n]^\Dimension)_{n\in\NN}$ are a van~Hove sequence.

\par
The \emph{Gilbert graph} of a configuration $C$ has vertices $C$ and edges connecting points at distance at most $\Radius$, i.e., $\DescSet{\EnumSet{x,y}\subseteq{}C}{\Distance{x-y}\le\Radius}$.
The configuration $C$ is a \emph{$\Radius$--cluster}, if it is \emph{$\Radius$--connected}, i.e., its Gilbert graph is connected.
Two points $x$ and $y$ are \emph{$\Radius$--connected by a configuration} $C$, written $x\RConnectedToIn{C}y$, if there is a finite path of jumps of at most  distance $\Radius$ between $x$ and $y$ using only points in $C$ as intermediate points.
Two Borel sets are $\Radius$--connected by a configuration $C$, if there is a $\Radius$--connected pair of points, with one point from each set.

\subsection{Point processes}
\label{sec_pp}

\par
For $B\in\Borel$, let $\Configs{B}$ be the \emph{locally finite point configurations} on $B$, i.e., for each $C\in\Configs{B}$ and $A\in\BoundedBorel$, $\Cardinality{C\cap{}A}<\infty$.
Let $\Algebra{B}$ be the $\sigma$--algebra on $\Configs{B}$ generated by $\DescSet{\DescSet{C\in\Configs{B}}{C\cap{}A=\emptyset}}{B\supseteq A\in\Borel}$, i.e., compatible with the \emph{Fell topology}.

\par
A \emph{simple point process} (short PP) on a Borel set $B\in\Borel$ is a random variable taking values in $\Configs{B}$.
This work treats a PP as a \emph{locally finite random subset of points} of $\RRD$, instead of as a random measure or as a collection of marginal counting rvs.
Let $\Law$ be a PP law and denote by $\xi$ the canonical variable on $\Configs{\RRD}$.
For $B\supseteq{}A\in\Borel$, abbreviate $\xi\cap{}A$ to $\xi_A$.

\par
A Borel measure $M$ on $(\Configs{B},\Algebra{B})$ is the \emph{local Janossy measure}~\cite[after (5.3.2)]{Daley_VereJones__AnIntroductionToTheTheoryOfPointProcesses_I__Springer_2003} of $\Law$ on $B\in\BoundedBorel$, if
\begin{equation}
\label{eq_janossy_measure}
 \forall\,E\in\Algebra{B}:
 \quad
 \Law(\xi_{B}\in{}E) = \int_E M(\D{C})
 \,.
\end{equation}
This definition of local Janossy measure is a portmanteau version of the traditional definition on generating cylinder sets.

\par
Because the local Janossy measure of a PP law $\Law$ on $B\in\Borel$ on $B\supseteq A\in\BoundedBorel$ equals the Janossy measure of the restriction of the law to $A$, the remainder of this paper drops the quantifier ``local''.
If $\xi$ has finite moment measures of all orders under $\Law$, then the Janossy measure in~\eqref{eq_janossy_measure} exists~\cite[Theorem 5.4.I]{Daley_VereJones__AnIntroductionToTheTheoryOfPointProcesses_I__Springer_2003}.
For $B\in\BoundedBorel$ and $C\in\Configs{B}$, write the infinitesimal of the Janossy measure of $\Law$ on $B$ at $C$ as $\Law(\xi_{B}=\D{C})$.

\par
The \emph{intensity measure} of the PP law $\Law$ is the average number of points on bounded Borel sets.
For $B\in\BoundedBorel$, it equals $\int_{\Configs{B}} \Cardinality{C} \Law(\xi_{B}=\D{C})$.

\subsection{The Boolean model}
\label{sec_boolean}

\par
The classic PP is the \emph{Poisson PP law} $\Poisson{B,\alpha}$ of intensity $\alpha$ on $B\in\Borel$, i.e., with intensity measure $\alpha\Lebesgue$.

\par
A configuration $C\in\Configs{\RRD}$ \emph{$\Radius$--percolates}, if it contains an infinite $\Radius$--cluster.
The bounded finiteness of $C$ renders this equivalent to the existence of an unbounded $\Radius$--cluster.
The \emph{Boolean model} of intensity $\alpha$ is a $\Poisson{\RRD,\alpha}$--distributed PP, with closed spheres of radius $\Radius/2$ centred at the points.
If spheres overlap, then the corresponding points are connected.
This is just $\Radius$--connectivity from Section~\ref{sec_space}.
The Boolean model \emph{percolates}, if it contains an infinite $\Radius$--cluster.

\par
Adding more points improves $\Radius$--connectivity.
Hence, the probability of percolation is monotone increasing in $\alpha$.
The Poissonian nature of the Boolean model makes percolation a tail event, i.e., it holds with either probability $0$ or $1$.
Thus, a critical intensity separates the non-percolating and percolating regimes.

\begin{theorem}[{\cite[Theorem 3.3]{Meester_Roy__ContinuumPercolation__CUP_1996}}]\label{thm_boolean}
For $\Dimension\ge 2$, a $\CritValuePercolation{\Dimension}\in{}]0,\infty[$ separates the sub-critical (almost-never percolating) from the super-critical (almost-surely percolating) intensities.
If $\alpha<\CritValuePercolation{\Dimension}$ and $(B_n)_{n\in\NN}$ is van~Hove, then
\begin{equation}
\label{eq_connection_subcritical}
 \Poisson{B_n,\alpha}(A\RConnectedToIn{\xi}\Ring{B_n})
 \xrightarrow[n\to\infty]{}0
 \,.
\end{equation}
\end{theorem}

\par
In one dimension, percolation almost-never happens at finite intensities.
Whence, $\CritValuePercolation{1}=\infty$~\cite[Theorem 3.1]{Meester_Roy__ContinuumPercolation__CUP_1996}.
In the sub-critical regime, the size of the $\Radius$--cluster containing the origin decays exponentially~\cite[Section 3.7]{Meester_Roy__ContinuumPercolation__CUP_1996}.
Section~\ref{sec_bounds} discusses bounds on $\CritValuePercolation{\Dimension}$.

\subsection{The hard-sphere model}
\label{sec_hs}

\par
Let $\Iverson{.}$ be \emph{Iverson brackets}\footnote{They work better with diverse logical expressions than indicator functions.}.
For disjoint $Y,C\in\Configs{\RRD}$, the indicator function $\SymbolIsHS$ of the \emph{conditional hard-core constraint} of $Y$ under condition $C$ is given by
\begin{equation}
\label{eq_ishs_def}
 \IsHS{Y|C}:=
 \prod_{\EnumSet{x,y}\subseteq{}Y} \Iverson{\Distance{x-y}>\Radius}
 \prod_{y\in{}Y,x\in{}C}
 \Iverson{\Distance{x-y}>\Radius}
 \,.
\end{equation}
For a \emph{bounded domain} $B\in\BoundedBorel$, a \emph{boundary condition} $C\in\Configs{\Complement{B}}$ and an \emph{activity} $\lambda\in[0,\infty[$, consider the \emph{hard-sphere model} with law $\HardSphere{B,C,\lambda}$.
As it is the Poisson PP of intensity $\lambda$ conditioned to be hard-core, its Janossy infinitesimal is
\begin{equation}
\label{eq_hsjan_def}
 \HardSphere{B,C,\lambda}(\D{Y})
 = \Poisson{B,\lambda}(\D{Y}|\IsHS{\xi|C}=1)
 \,.
\end{equation}
\begin{subequations}
\label{eq_statmech}
The alternative definition in statistical mechanics uses the \emph{pair potential}
\begin{equation}
\label{eq_statmech_potential}
 \Potential:\quad(\RR^\Dimension)^2\mapsto[0,\infty]
 \quad
 (x,y)\mapsto
 \begin{cases}
 \infty
 &\text{if }\Distance{x-y}\le\Radius\,,
 \\
 0
 &\text{if }\Distance{x-y}>\Radius\,.
 \end{cases}
\end{equation}
The Hamiltonian of $n$ ordered points in $B$ is
\begin{equation}
\label{eq_statmech_hamiltonian}
 H(x_1,\dotsc,x_n|C)
 :=
 \sum_{1\le{}i<j\le{}n} \Potential(x_i,x_j)
 +
 \sum_{1\le{}i\le n, y\in{}C} \Potential(x_i,y)
 \,.
\end{equation}
The density of $x\in{}B^n$ is
\begin{equation}
\label{eq_statmech_density}
 \HardSphere{B,C,\lambda}(\D{x})
 := \frac{\lambda^n e^{-H(x|C)}}{n! \PartFun{B,C,\lambda}}  \D{x}
 \,,
\end{equation}
where the \emph{partition function} $\SymbolPartFun$ is
\begin{equation}
\label{eq_statmech_partfun}
 \PartFun{B,C,\lambda}
 := \sum_{n=0}^\infty \frac{\lambda^n}{n!}
  \int_{B^n}  e^{-H(x|C)} \D{x}
 \,.
\end{equation}
The convention $e^{-\infty}=0$ encodes~\eqref{eq_ishs_def} by~\eqref{eq_statmech_hamiltonian}.
The remainder of this paper uses the PP notation as in~\eqref{eq_hsjan_def}, except for the partition function $\SymbolPartFun$.
\end{subequations}
Because of the bounded range interaction in $\IsHS{Y|C}$ in~\eqref{eq_ishs_def}, one may restrict the boundary condition to $\Configs{\Ring{B}}$.

\par
A \emph{Gibbs measure} is a weak limit of a sequence $(\HardSphere{B_n,C_n,\lambda})_{n\in\NN}$ along a van~Hove sequence $(B_n)_{n\in\NN}$ and a sequence $(C_n)_{n\in\NN}$ of boundary conditions with $C_n\in\Configs{\Complement{B_n}}$~\cite[Sections 2 and 3]{Preston__RandomFields__LNM_Springer_1976}.
The Gibbs measures $\GibbsMeasures{\lambda}$ of the \emph{specification} $\HardSphere{\lambda}:=(\HardSphere{B,C,\lambda})_{B\in\BoundedBorel,C\in\Configs{\Complement{B}}}$ form a simplex.
Unlike in the lattice case~\cite{Runnels__PhaseTransitionsOfHardSphereLatticeGases__CMP_1975}, in the continuum case of dimension greater than one, the existence of a finite critical activity at which a phase transition happens is widely believed, but not yet proven.
See the solution in one dimension~\cite{Tonks__TheCompleteEquationOfStateOfOneTwoAndThreeDimensionalGasesOfElasticHardSpheres__PhysRev_1936}, the absence of positional phase transition in two dimensions~\cite{Richthammer__TranslationInvarianceOfTwoDimensionalGibbsianPointProcesses__CMP_2007}, which does not exclude a conjectured orientational phase transition, and the discussion of the state of the problem in higher dimensions~\cite[Section 3.3]{Loewen__FunWithHardSpheres__LNP_1999}.
If $\Radius=0$, then there is no interaction, the hard-sphere model reduces to a Poisson PP and $\Poisson{\RRD,\lambda}$ is the unique Gibbs measure.

\subsection{Stochastic domination}
\label{sec_stoch_dom}

\par
On $\Configs{B}^n$, the standard product $\sigma$--algebra is $\ProductAlgebra{B}{n}$.
The canonical variables on $\Configs{B}^n$ are $\xi:=(\xi^1,\dotsc,\xi^n)$.
A \emph{coupling} $\Law$ of $n$ PP laws $\Law_1,\dotsc,\Law_n$ on $B\in\Borel$ is a probability measure on $(\Configs{B}^n,\ProductAlgebra{B}{n})$ such that, for all $1\le{}i\le{}n$ and $E\in\Algebra{B}$, $\Law(\xi^i\in{}E)=\Law_i(\xi\in{}E)$.

\par
A PP law $\Law_2$ \emph{stochastically dominates} a PP $\Law_1$, if there exists a coupling $\Law$ of them with $\Law(\xi^1\subseteq\xi^2)=1$.
Equivalently, by Strassen's theorem, for all increasing events $E$, $\Law_1(E)\le\Law_2(E)$.
A Poisson PP stochastically dominates a hard-sphere model with the same activity as the intensity of the Poisson PP~\cite[Example 2.2]{Georgii_Kueneth__StochasticComparisonOfPointRandomFields__JAP_1997}.

\section{Results}
\label{sec_results}

\subsection{Disagreement percolation}
\label{sec_daperc}

\par
At the core of disagreement percolation is a coupling of two instances of the hard-sphere model on the same finite volume, but with differing boundary conditions, such that the set of points differing between the two instances (the \emph{disagreement cluster}) is stochastically dominated by a Poisson PP.
Therefore, one may control the disagreement clusters and the influence of the differing boundary conditions by the percolation clusters of the Boolean model.

\par
If the intensity of the dominating Poisson PP is below the critical value for percolation in the Boolean model, then the finiteness of percolation clusters controls the influence of the differing boundary conditions.
The influence vanishes as the finite volume tends to the whole space.
This implies the uniqueness of the Gibbs measure of the hard-sphere model.
Furthermore, as the cluster size of the Boolean model decays exponentially in the sub-critical phase, controls of the Gibbs measure such as the influence of boundary conditions or the reduced pair correlation function decay exponentially, too.

\par
The remainder of this section formalises the preceding outline.
The proofs are in Section~\ref{sec_daperc_proof}.
The \emph{symmetric difference} $S_1\SymDiff{}S_2$ between sets $S_1$ and $S_2$ equals $(S_1\setminus{}S_2)\cup(S_2\setminus{}S_1)$.

\begin{definition}
\label{def_dac}
Let $\alpha,\lambda\in[0,\infty[$.
A \emph{disagreement coupling} on $B\in\BoundedBorel$ with $C_1,C_2\in\Configs{\Complement{B}}$ of intensity $\alpha$ and activity $\lambda$ is a law $\Law$ on $(\Configs{B}^3,\ProductAlgebra{B}{3})$ with
\begin{subequations}
\label{eq_dac_def}
\begin{gather}
\label{eq_dac_def_hs}
  \forall\,1\le{}i\le{}2, E\in\Algebra{B}:\quad
  \Law(\xi^i\in{}E)=\HardSphere{B,C_i,\lambda}(\xi\in{}E)
 \,,\\\label{eq_dac_def_poi}
  \forall\,E\in\Algebra{B}:\quad
  \Law(\xi^3\in{}E)=\Poisson{B,\alpha}(\xi\in{}E)
 \,,\\\label{eq_dac_def_subset}
  \Law(\xi^1\SymDiff\xi^2\subseteq\xi^3)=1
 \,,\\\label{eq_dac_def_connected}
  \Law(
    \forall\,x\in\xi^1\SymDiff\xi^2:
    x\RConnectedToIn{\xi^1\SymDiff \xi^2}C_1\SymDiff C_2
   )=1
 \,.
\end{gather}
\end{subequations}
A \emph{disagreement coupling family} of intensity $\alpha$ and activity $\lambda$ is a family of disagreement couplings $(\Law_{B,C_1,C_2,\lambda,\alpha})_{B\in\BoundedBorel,C_1,C_2\in\Configs{\Complement{B}}}$.
\end{definition}

\par
A disagreement coupling family in the sub-critical phase of the Boolean model implies uniqueness of the Gibbs measure.

\begin{theorem}
\label{thm_uniqueness}
If there exists a disagreement coupling family of intensity $\alpha<\CritValuePercolation{\Dimension}$ at activity $\lambda$, then $\GibbsMeasures{\lambda}$ consists of a single Gibbs measure.
\end{theorem}

\par
Disagreement percolation also implies that the \emph{sensitivity to changes in the boundary condition}~\eqref{eq_exp_decay_sbc}, the \emph{finite volume error}~\eqref{eq_exp_decay_fve} and the \emph{probabilities of separated events}~\eqref{eq_exp_decay_sepev} on small sets decay exponentially.
The rate of exponential decay is the same as the one of the Boolean model~\eqref{eq_exp_decay_boolean}, which holds in the whole sub-critical regime of the Boolean model~\cite[Section 3.7]{Meester_Roy__ContinuumPercolation__CUP_1996}.

\begin{theorem}
\label{thm_exp_decay}
Assume the existence of a disagreement coupling family of intensity $\alpha<\CritValuePercolation{\Dimension}$ at activity $\lambda$ and
a there exist $K\ge{}1,\kappa>0$ such that, for all $A,B\in\BoundedBorel$ with $\Diameter{A}\le{}1$,
\begin{equation}
\label{eq_exp_decay_boolean}
 \Poisson{\RRD,\alpha}(A\RConnectedToIn{\xi} B)
 \le{}K e^{-\kappa{}\SetDistance{A,B}}
 \,.
\end{equation}
\begin{subequations}
\label{eq_exp_decay_controls}
For all $A,B\in\BoundedBorel$ with $\Diameter{A}\le{}1$, $A\subseteq{}B$, $C\in\Configs{\Complement{B}}$, $x\in\Complement{(B\cup{}C)}$ and $E\in\Configs{A}$,
\begin{equation}
\label{eq_exp_decay_sbc}
 \Modulus{\HardSphere{B,C,\lambda}(\xi_{A}\in{}E)
  -\HardSphere{B,C\cup\EnumSet{x},\lambda}(\xi_{A}\in{}E)
 }
 \le K e^{-\kappa{}\SetDistance{A,\EnumSet{x}}}
 \,.
\end{equation}
Let $\nu$ be the unique Gibbs measure in $\GibbsMeasures{\lambda}$.
For all $A,B\in\BoundedBorel$ with $\Diameter{A}\le{}1$, $A\subseteq{}B$, $C\in\Configs{\Complement{B}}$ and $E\in\Configs{A}$,
\begin{equation}
\label{eq_exp_decay_fve}
 \Modulus{\HardSphere{B,C,\lambda}(\xi_{A}\in{}E)-\nu(\xi_{A}\in{}E)}
 \le K e^{-\kappa{}\SetDistance{A,\Complement{B}}}
 \,.
\end{equation}
For all $A,B\in\BoundedBorel$ with $\Diameter{A}\le{}1$, $E\in\Configs{A}$ and $F\in\Configs{B}$,
\begin{equation}
\label{eq_exp_decay_sepev}
 \Modulus{\nu(\xi_{A}\in{}E,\xi_{B}\in{}F)
  -\nu(\xi_{A}\in{}E)\nu(\xi_{B}\in{}F)
 }
 \le K e^{-\kappa{}\SetDistance{A,B}}
 \,.
\end{equation}
\end{subequations}
\end{theorem}

\subsection{Bounds from disagreement percolation}
\label{sec_bounds}

The hard-sphere model admits a disagreement coupling family of the same intensity as its activity.

\begin{theorem}
\label{thm_twisted}
There exists a disagreement coupling family of intensity $\lambda$ for $\HardSphere{\lambda}$, with property~\eqref{eq_dac_def_subset} improved to
\begin{equation}
\label{eq_dac_twisted_subset}
 \Law(\xi^1\cup{}\xi^2\subseteq\xi^3)=1
 \,.
\end{equation}
If $\lambda<\CritValuePercolation{\Dimension}$, then $\GibbsMeasures{\lambda}$ is a singleton and exponential decay as in~\eqref{eq_exp_decay_controls} holds.
\end{theorem}

Theorem~\ref{thm_twisted} follows from the disagreement coupling family in Section~\ref{sec_twisted} and theorems~\ref{thm_uniqueness} and~\ref{thm_exp_decay}.
A motivation of this coupling is in Section~\ref{sec_motivation} and discussion of generalisations and other approaches in Section~\ref{sec_outlook}.

\begin{subequations}
\label{eq_critin}
\par
Bounds on $\CritValuePercolation{\Dimension}$ translate directly into sufficient conditions for the uniqueness of the Gibbs measure.
In one dimension the Boolean model never percolates~\cite[Theorem 3.1]{Meester_Roy__ContinuumPercolation__CUP_1996}.
\begin{equation}
\label{eq_critin_1d}
 \CritValuePercolation{1}=\infty
 \,.
\end{equation}
In two dimensions, rigorous bounds on $\CritValuePercolation{2}$ are $[\frac{0.174}{\Radius^2},\frac{0.843}{\Radius^2}]$~\cite[Theorem 3.10]{Meester_Roy__ContinuumPercolation__CUP_1996}.
More recent high confidence bounds in~\cite{Balister_Bollobas_Walters__ContinuumPercolationWithStepsInTheSquareOrTheDisc__RSA_2005}, taken from~\cite[Equation (2)]{Mertens_Moore__ContinuumPercolationThresholdsInTwoDimensions__PhysRevE_2012}, are
\begin{equation}
\label{eq_critin_2d}
 \frac{0.358}{\Radius^2}
 \sim
 \frac{1.127}{\pi \Radius^2}
 <
 \CritValuePercolation{2}
 \,.
\end{equation}
For dimensions $2$ to $10$, simulation bounds are in ~\cite{Jiao_Torquato__EffectOfDimensionalityOnTheContinuumPercolationOfOverlappingHyperspheresAndHypercubes_II_SimulationResultsAndAnalysis__JChemPhys_2012,Jiao_Torquato__Erratum_EffectOfDimensionalityOnTheContinuumPercolationOfOverlappingHyperspheresAndHypercubes_II_SimulationResultsAndAnalysis__JChemPhys_2014}.
Another set of high confidence and rigorous bounds via an Ornstein-Zernike approach are in~\cite[Table 4]{Ziesche__SharpnessOfThePhaseTransitionAndLowerBoundsForTheCriticalIntensityInContinuumPercolationOnRD__AIHP_2018}.
The asymptotic behaviour of the critical intensity~\cite{Penrose__ContinuumPercolationAndEuclideanMinimalSpanningTreesInHighDimensions__AAP_1996}, taken from~\cite[Section 3.10]{Meester_Roy__ContinuumPercolation__CUP_1996}, is
\begin{equation}
\label{eq_critin_asymptotic}
 \lim_{\Dimension\to\infty}
  \CritValuePercolation{\Dimension}
  \VolumeCoeff{\Dimension}\Radius^\Dimension
 = 1
 \,.
\end{equation}
\end{subequations}

The inequality~\eqref{eq_exp_decay_sepev} from Theorem~\ref{thm_twisted} implies the exponential decay of the \emph{reduced pair correlation function}.

\begin{theorem}
\label{thm_exp_decay_rpcf}
If $\lambda<\CritValuePercolation{\Dimension}$ and~\eqref{eq_exp_decay_boolean} holds with $\alpha=\lambda$, then, for all $x,y\in\RRD$ and using the constants from~\eqref{eq_exp_decay_boolean}, the reduced pair correlation function $\ReducedPairCorrelation$ decays as
\begin{equation}
\label{eq_exp_decay_rpcf}
 \ReducedPairCorrelation(x,y)
 \le
 K e^{-\kappa{}\Distance{x-y}}
 \,.
\end{equation}
\end{theorem}

\subsection{Comparison with expansion bounds}
\label{sec_comparison}

\par
Popular methods to study the absence of phase transitions, in particular to guarantee the uniqueness of the Gibbs measure, are virial and cluster expansion methods~\cite{Ruelle__StatisticalMechanics_RigorousResults__Ben_1969}.
Both deliver analyticity of the free energy, too.
Let $\CERadius{\Dimension}$ be the radius of the cluster expansion in $\Dimension$ dimensions.

\begin{subequations}
\label{eq_ce}
\par
In one dimension, disagreement percolation~\eqref{eq_critin_1d} replicates Tonks' classic result of the complete absence of phase transitions via virial expansion methods~\cite{Brydges_Imbrie__DimensionalReductionFormulasForBranchedPolymerCorrelationFunctions__JSP_2003,Jansen__ClusterAndVirialExpansionsForTheMultiSpeciesTonksGas__JSP_2015,Labelle_Leroux_Ducharme__GraphWeightsArisingFromMayersTheoryOfClusterIntegrals__SemLotCom_2005,Tonks__TheCompleteEquationOfStateOfOneTwoAndThreeDimensionalGasesOfElasticHardSpheres__PhysRev_1936}.
In terms of the activity, it is known that the radius of the cluster expansion is exactly~\cite{Brydges_Imbrie__DimensionalReductionFormulasForBranchedPolymerCorrelationFunctions__JSP_2003,HoferTemmel__ShearersPointProcessAndTheHardSphereModelInOneDimension,Labelle_Leroux_Ducharme__GraphWeightsArisingFromMayersTheoryOfClusterIntegrals__SemLotCom_2005}
\begin{equation}
\label{eq_ce_1d}
 \CERadius{1}=\frac{1}{e\Radius}
 \,.
\end{equation}

\par
In two dimensions, the best currently known lower bounds~\cite{Fernandez_Procacci_Scoppola__TheAnalyticityRegionOfTheHardSphereGas_ImprovedBounds__JSP_2007} and upper bounds~\cite[Section 4.5]{Ruelle__StatisticalMechanics_RigorousResults__Ben_1969}
are
\begin{equation}
\label{eq_ce_2d}
 \frac{0.1625}{\Radius^2}
 \sim
 \frac{0.5107}{\pi\Radius^2}
 <
 \CERadius{2}
 <
 \frac{2}{e\pi\Radius^2}
 \sim
 \frac{0.2342}{\Radius^2}
 \,.
\end{equation}
The bounds in~\eqref{eq_ce_2d} are between $0.45$ and $0.65$ times the disagreement percolation bound in~\eqref{eq_critin_2d}.
General bounds on the cluster expansion radius from~\cite[Section 4.5]{Ruelle__StatisticalMechanics_RigorousResults__Ben_1969} are
\begin{equation}
\label{eq_ce_general}
 \frac{1}{e\VolumeCoeff{\Dimension}\Radius^\Dimension}
 \le
 \CERadius{\Dimension}
 \le
 \frac{2}{\VolumeCoeff{\Dimension}\Radius^\Dimension}
 \,.
\end{equation}
As $\VolumeCoeff{1}=2$, equation~\eqref{eq_ce_1d} shows that the upper bound is tight.
I conjecture that the asymptotic behaviour in high dimensions is
\begin{equation}
\label{eq_ce_asymptotic}
 \lim_{\Dimension\to\infty}
  \CERadius{\Dimension}\VolumeCoeff{\Dimension}\Radius^\Dimension
 = \frac{1}{e}\,.
\end{equation}
\end{subequations}
Comparing~\eqref{eq_critin_asymptotic} and~\eqref{eq_ce_asymptotic}, the asymptotic improvement should be by a factor of $e$.
This is not surprising, because on the infinite $k$--regular tree $T_k$, the critical percolation probability is $\frac{1}{k-1}$ and the radius of the cluster expansion is $\frac{(k-2)^{k-2}}{(k-1)^{k-1}}\sim\frac{1}{e(k-1)}$.
Both $\ZZ^\Dimension$ and $\RR^\Dimension$ behave for large $\Dimension$ as $T_{2\Dimension}$, for both percolation and cluster expansion.
Extrapolating arguments of~\cite[Section 8]{Scott_Sokal__TheRepulsiveLatticeGasTheIndependentSetPolynomialAndTheLovaszLocalLemma__JSP_2005} gives a heuristic for the upper bound in~\eqref{eq_ce_general}, too.
Finally, I conjecture that disagreement percolation is always better than cluster expansion.
A possible approach is recent work connecting the Ornstein-Zernike equation for the Boolean model with Ruelle-like sufficient conditions for cluster expansion~\cite{Last_Zieschke__OnTheOrnsteinZernikeEquationForStationaryClusterProcessesAndTheRandomConnectionModel__AdvAppProb_2017}.

\subsection{Motivation behind the dependent thinning and twisted coupling}
\label{sec_motivation}

\par
This section assumes familiarity with the dependent thinning in Definition~\ref{def_thin} and the twisted disagreement coupling family in Definition~\ref{def_twisted}.

The approach to disagreement percolation in~\cite{VanDenBerg_Maes__DisagreementPercolationInTheStudyOfMarkovFields__AP_1994} is a vertex-wise conditional coupling of two Markov fields on a finite graph.
A uniform control of those couplings allows stochastic domination by a Bernoulli product field.
This poses a problem on $\RRD$.
The key insight is to flip the picture around.
Start with the Bernoulli random field and reinterpret the conditional couplings as simultaneous dependent thinnings to the two dominated Markov fields.
Transferring this to the PP case is non-trivial, but helpfully~\cite{VanDenBerg_Steif__PercolationAndTheHardCoreLatticeGasModel__SPA_1994} introduced an optimisation for the hard-core model.
This reduces the thinning probability onto two hard-core models on a single vertex to a thinning probability of a single hard-core model.
The equivalent of the discrete thinning probability is the rhs of~\eqref{eq_derivative}.
This was the starting point of the generalisation to the PP case.
In the PP case, this enables the independent construction on the disjoint domains in~\eqref{eq_twisted_def_zone}.
It allows to ``twist'' two hard-sphere models of activity $\lambda$ under a single $\Poisson{\lambda}$ PP, i.e., have joint stochastic domination in~\eqref{eq_twisted_zone_dominate}.

\par
The overall recursive approach from the dependent vertex-wise couplings stays and translates into the recursive definition~\eqref{eq_twisted_def_rec}.
The recursive definition of $\TwistedRec{B,C_1,C_2}$ demands that it is jointly measurable in the boundary conditions $C_1$ and $C_2$.
By the above outline, this comes back to the measurability of $\HSPoiThin{B,C}$ in the boundary condition $C$.
The classic dominating couplings between a Poisson PP and a hard-sphere model in~\cite{Preston__SpatialBirthAndDeathProcesses__BIIS_1975, Georgii_Kueneth__StochasticComparisonOfPointRandomFields__JAP_1997} are implicit and do not provide this measurability readily.
Besides the generalisation of the thinning probability from the discrete case, this is the main reason for the explicit construction of $\HSPoiThin{B,C}$ in Section~\ref{sec_thin_def}.
But the $D=\emptyset$ case in~\eqref{eq_twisted_def_rec} suggests to use the dependent thinning approach for a single dominated hard-sphere model, too.
In this case, the calculations are doable and lead to the dependent thinning in Definition~\ref{def_thin}.

\par
The \emph{Papangelou intensity}~\cite[(15.6.13)]{Daley_VereJones__AnIntroductionToTheTheoryOfPointProcesses_II__Springer_2008} is the infinitesimal cost of adding another point to a given configuration.
It is $\IsHS{\EnumSet{x}|Y\cup{}C}\lambda$ for the hard-sphere model.
The Poisson PP has constant Papangelou intensity $\lambda$.
Thus, one can control the hard-sphere model pointwise incrementally by a Poisson PP.
All three stochastic dominations of a hard-sphere model by a Poisson PP (the dependent thinning in Definition~\ref{def_thin},~\cite{Georgii_Kueneth__StochasticComparisonOfPointRandomFields__JAP_1997} and~\cite{Preston__SpatialBirthAndDeathProcesses__BIIS_1975}) build upon this fact.
In the $D=\emptyset$ case, $\TwistedRec{}$ reduces to the same setting, too.
It is yet unknown if this is the smallest Poisson intensity needed to dominate the hard-sphere model.
Looking at~\eqref{eq_thin_proba} and~\eqref{eq_derivative}, one sees that $\lambda$ is indeed approached for $x$ far enough away from $C$ and $Y$ and large $B$.
The rewrite of the thinning probability as the derivative of a finite volume free energy in~\eqref{eq_derivative} adds another interpretation: acceptance of a point $x$ happens with a probability depending on the change in the free energy.
For large domains the point would change nothing in the free energy, whence it could be accepted with probabilities approaching $1$.
Again, no smaller Poisson intensity smaller than $\lambda$ allows this.
Because the thinning procedure depends on the ordering from Section~\ref{sec_ordering_derivative}, I consider the preceding thoughts only a strong indicator but not a proof for the minimality of $\lambda$.

\par
Another natural question is whether the depending thinning factorises over clusters of the dominating Poisson PP.
Although it looks likely to be true, because the answer is not relevant here, this question is not investigated.

\subsection{Outlook}
\label{sec_outlook}

\par
In the lattice case, disagreement percolation implies the complete analyticity of the free energy, pointed out by Schonmann~\cite[Note added in proof]{VanDenBerg_Maes__DisagreementPercolationInTheStudyOfMarkovFields__AP_1994}, and the Poincare inequality for the usual spin-flip dynamics~\cite{Chazottes_Redig_Voellering__ThePoincareInequalityForMarkovRandomFieldsProvedViaDisagreementPercolation__IMNS_KNAW_2011}.
In principle, both results should be generalisable to the hard-sphere model, too.
The exponential control in~\eqref{eq_exp_decay_controls} looks exactly like what is needed in the discrete case for complete analyticity~\cite{Dobrushin_Shlosman__CompletelyAnalyticalInteractions_AConstructiveDescription__JStatPhys_1987}, but a theory for PPs is still missing.

\par
The proofs of theorems~\ref{thm_uniqueness} and~\ref{thm_exp_decay} are independent of the hard-sphere model and apply to arbitrary Gibbs PPs with bounded range interaction.
A generalisation to the physically interesting case of marked Gibbs PP models with finite, but unbounded, range should be possible.
This demands a notational and definitional base exceeding the limits of a single paper, though.

\par
Beyond the hard-sphere model, one could do a product construction in~\eqref{eq_twisted_def_zone} and compensate by adding an additional $\Poisson{B,\lambda}$ in the $D=\emptyset$ case in~\eqref{eq_twisted_def_zone}.
This would lead to a disagreement coupling family of intensity $2\lambda$, for a repulsive potential.
The recursive construction still demands the dominating coupling to be measurable in the boundary conditions.

\par
The more simple product approach from~\cite{VanDenBerg__AUniquenessConditionForGibbsMeasuresWithApplictionToThe2DimensionalIsingAntiferromagnet__CMP_1993} with a swapping argument yields only a lower bound of $\CritValuePercolation{\Dimension}/2$.
Thus, it is not strong enough for the comparison in Section~\ref{sec_comparison}.
Also, the same measurability concerns as in the twisted approach surface, too.

\par
Another sufficient condition for uniqueness Gibbs measure, and even complete analyticity of the free energy, is Dobrushin's uniqueness condition~\cite{Dobrushin_Shlosman__ConstructiveCriterionForTheUniquenessOfGibbsField__SPDS_1985}.
There have been generalisations to the PP case~\cite{Klein__DobrushinUniquenessTechniquesAndTheDecayOfCorrelationsInContinuumStatisticalMechanics__CMP_1982,Klein__ConvergenceOfGrandCanonicalGibbsMeasures__CMP_1984}, but I make no explicit comparison here.

\section{Proof of theorems~\ref{thm_uniqueness},~\ref{thm_exp_decay} and~\ref{thm_exp_decay_rpcf}}
\label{sec_daperc_proof}

\par
The proof of Theorem~\ref{thm_uniqueness} follows closely the one in the discrete case~\cite[proof of corollaries 1 and 2]{VanDenBerg_Maes__DisagreementPercolationInTheStudyOfMarkovFields__AP_1994}.
Proposition~\ref{prop_diff_connection_bound} applies a disagreement coupling to bound the difference between the two hard-sphere models by a percolation connection probability.
This proposition is the key control of the influence of the differing boundary conditions.
Theorem~\ref{thm_uniqueness} uses a disagreement coupling family to exploit these bounds on increasing scales.
First, it restricts to a small domain, then it applies the bounds from disagreement coupling and finally, it uses the sub-criticality of the Boolean model to tighten the bound to zero as the domain increases.
Theorem~\ref{thm_exp_decay} uses Proposition~\ref{prop_diff_connection_bound} to control the influence of the differing boundary conditions on general events.
Theorem~\ref{thm_exp_decay_rpcf} derives a tighter disagreement bound for increasing events from~\eqref{eq_dac_twisted_subset} and proves exponential decay of the pair correlation function.

\begin{proposition}
\label{prop_diff_connection_bound}
Let $A,B\in\BoundedBorel$ with $A\subseteq{}B$, $C_1,C_2\in\Configs{\Complement{B}}$ and $\alpha,\lambda\in[0,\infty[$.
Let $\Law:=\Law_{B,C_1,C_2,\lambda,\alpha}$ be a disagreement coupling.
For $E\in\Algebra{A}$,
\begin{equation}
\label{eq_diff_connection_bound}
 \Modulus{
  \HardSphere{B,C_1,\lambda}(\xi_{A}\in{}E)
  -\HardSphere{B,C_2,\lambda}(\xi_{A}\in{}E)
 }
 \le
 \Poisson{B,\alpha}(A\RConnectedToIn{\xi}C_1\SymDiff{}C_2)
 \,.
\end{equation}
\end{proposition}

\begin{proof}
First, reduce the difference by cancelling symmetric parts.
\begin{align*}
 &\Modulus{
  \HardSphere{B,C_1,\lambda}(\xi_{A}\in{}E)
  -\HardSphere{B,C_2,\lambda}(\xi_{A}\in{}E)
  }
 \\\stackrel{\eqref{eq_dac_def_hs}}{=}{}
 &\Modulus{\Law(\xi^1_{A}\in{}E)-\Law(\xi^2_{A}\in{}E)}
 \\={}
 &\Modulus{\Law(\xi^1_{A}\in{}E,\xi^2_{A}\not\in{}E)
  -\Law(\xi^1_{A}\not\in{}E,\xi^2_{A}\in{}E)}
 \\\le{}
 &\max\EnumSet{
  \Law(\xi^1_{A}\in{}E,\xi^2_{A}\not\in{}E)
  ,\Law(\xi^1_{A}\not\in{}E,\xi^2_{A}\in{}E)
  }
 \,.
\end{align*}
Second, relax the asymmetry to disagreement and use disagreement percolation.
\begin{align*}
 \Law(\xi^1_{A}\in{}E,\xi^2_{A}\not\in{}E)
 \stackrel{\text{relax}}{\le{}}
 &\Law(\xi^1_{A}\SymDiff\xi^2_{A}\not=\emptyset)
 \\\stackrel{\eqref{eq_dac_def_connected}}{\le}{}
 &\Law(
  A\RConnectedToIn{\xi^1\SymDiff\xi^2}C_1\SymDiff C_2
  )
 \\\stackrel{\eqref{eq_dac_def_subset}}{\le}{}
 &\Law(
  A\RConnectedToIn{\xi^3}C_1\SymDiff C_2
  )
 \\\stackrel{\eqref{eq_dac_def_poi}}{=}{}
 &\Poisson{B,\alpha}(A\RConnectedToIn{\xi}C_1\SymDiff C_2)
 \,.\qedhere
\end{align*}
\end{proof}

\par
For disjoint $A,B\in\BoundedBorel$ and $E\in\Algebra{B}$, the following identities hold for the Janossy infinitesimals.
\begin{equation}
\label{eq_janossy_identities}
\begin{aligned}
 \Law(\xi_{A}=\D{Y})={}
 &
 \int_{\Configs{B}} \Law(\xi_{A\cup{}B}=\D{(Y\cup{}Z)})
 \,,\\
 \Law(\xi_{A}=\D{Y},\xi_{B}\in{}E)={}
 &
 \int_{\Configs{B}}\Iverson{Z\in{}E} \Law(\xi_{A\cup{}B}=\D{(Y\cup{}Z)})
 \,,\\
 \Law(\xi_{A}=\D{Y}|\xi_{B}\in{}E)={}
 &
 \int_{\Configs{B}}
   \frac{\Iverson{Z\in{}E}\Law(\xi_{A\cup{}B}=\D{(Y\cup{}Z)})}
        {\Law(\xi_{B}\in{}E)}
 \,.
\end{aligned}
\end{equation}

\begin{proof}[Proof of Theorem~\ref{thm_uniqueness}]
Let $\nu_1,\nu_2\in\GibbsMeasures{\lambda}$.
Showing that $\nu_1=\nu_2$ is equivalent to
\begin{equation*}
 \forall\,A\in\BoundedBorel, E\in\Algebra{A}:
 \quad
 \nu_1(\xi_{A}\in{}E)=\nu_2(\xi_{A}\in{}E)
 \,.
\end{equation*}
The hard-sphere property ensures that a Gibbs measure in $\GibbsMeasures{\lambda}$ has moment measures of all orders~\cite[(5.4.9)]{Daley_VereJones__AnIntroductionToTheTheoryOfPointProcesses_I__Springer_2003}.
Thus, its local Janossy measures exist.

\par
The following result controls the difference between two measures.
Let $\mu_1$ and $\mu_2$ be probability measures on the measurable space $(\Omega,\mathcal{A})$. For all $f:\Omega\to[0,1]$ measurable,
\begin{equation}
\label{eq_intsupdiff}
\begin{aligned}
 \LRModulus{\int f \D{\mu_1} - \int f \D{\mu_2} }
 ={}
 &\LRModulus{
  \int_{\Omega^2} (f(\omega_1)-f(\omega_2))
   \D{\mu_1}(\omega_1)\D{\mu_2}(\omega_2)
 }
 \\\le{}
 &\int_{\Omega^2}
   \Modulus{f(\omega_1)-f(\omega_2)}
   \D{\mu_1}(\omega_1)\D{\mu_2}(\omega_2)
 \\\le{}
 &\sup\DescSet{\Modulus{f(\omega_1)-f(\omega_2)}}{\omega_1,\omega_2\in\Omega}
 \,.
\end{aligned}
\end{equation}

\par
Let $(B_n)_{n\in\NN}$ be a van~Hove sequence with $A\subseteq{}B_1$.
For each Gibbs measure $\nu\in\GibbsMeasures{\lambda}$ and $n\in\NN$, the Gibbs property restricts the discussion to the bounded Borel set $B_n$.
Second, the existence of a disagreement coupling family of intensity $\alpha$ and~\eqref{eq_intsupdiff} controls the difference between different Gibbs measures by the connection probability of the Boolean model.
Taking the limit along the van~Hove sequence shows that the difference is zero.

\begin{align*}
 &\Modulus{\nu_1(\xi_{A}\in{}E)-\nu_2(\xi_{A}\in{}E)}
 \\\stackrel{\eqref{eq_janossy_identities}}{=}{}
 &\left|
    \int_{\Configs{\Complement{B_n}}}
    \HardSphere{B_n,C_1,\lambda}(\xi_{A}\in{}E)
    \nu_1(\xi_{\Complement{B_n}}=\D{C_1})
  \right.
 \\&-\left.
     \int_{\Configs{\Complement{B_n}}}
     \HardSphere{B_n,C_2,\lambda}(\xi_{A}\in{}E)
     \nu_2(\xi_{\Complement{B_n}}=\D{C_2})
   \right|
 \\\stackrel{\eqref{eq_intsupdiff}}{\le}{}
 &\sup\LRDescSet%
   {\LRModulus{\HardSphere{B_n,C_1,\lambda}(\xi_{A}\in{}E)
   -\HardSphere{B_n,C_2,\lambda}(\xi_{A}\in{}E)}}
   {C_1,C_2\in\Configs{\Complement{B_n}}}
 \\\stackrel{\eqref{eq_diff_connection_bound}}{\le}{}
 &\sup\LRDescSet%
   {\Poisson{B_n,\alpha}(A\RConnectedToIn{\xi}C_1\SymDiff C_2)}
   {C_1,C_2\in\Configs{\Complement{B_n}}}
 \\\stackrel{\text{relax}}{\le}{}
 &\Poisson{B_n,\alpha}(A\RConnectedToIn{\xi}\Complement{B_n})
 \xrightarrow[n\to\infty]{\eqref{eq_connection_subcritical}}{}
 0
 \,.\qedhere
\end{align*}
\end{proof}

\begin{proof}[Proof of Theorem~\ref{thm_exp_decay}]
For~\eqref{eq_exp_decay_sbc}, let $C_1:=C$ and $C_2:=C\cup\EnumSet{x}$.
Thus,
\begin{equation*}
 \Modulus{
  \HardSphere{B,C,\lambda}(\xi_{A}\in{}E)
  -\HardSphere{B,C\cup\EnumSet{x},\lambda}(\xi_{A}\in{}E)
 }
 \stackrel{\eqref{eq_diff_connection_bound}}{\le}{}
 \Poisson{B,\alpha}(A\RConnectedToIn{\xi}\EnumSet{x})
 \stackrel{\eqref{eq_exp_decay_boolean}}{\le}{}
 K e^{-\kappa{}\SetDistance{A,\EnumSet{x}}}
 \,.
\end{equation*}
For~\eqref{eq_exp_decay_fve}, one has
\begin{align*}
 &\Modulus{\HardSphere{B,C,\lambda}(\xi_{A}\in{}E)-\nu(\xi_{A}\in{}E)}
 \\\stackrel{\eqref{eq_janossy_identities}}{=}{}
 &\LRModulus{
  \HardSphere{B,C,\lambda}(\xi_{A}\in{}E)
  -\int_{\Configs{\Ring{B}}}
   \HardSphere{B,C',\lambda}(\xi_{A}\in{}E)
   \nu(\xi_{\Ring{B}}=\D{C'})
 }
 \\\le{}
 &\int_{\Configs{\Ring{B}}}
   \Bigl|
    \HardSphere{B,C,\lambda}(\xi_{A}\in{}E)
    -\HardSphere{B,C',\lambda}(\xi_{A}\in{}E)
   \Bigr|
   \nu(\xi_{\Ring{B}}=\D{C'})
 \\\stackrel{\eqref{eq_diff_connection_bound}}{\le}{}
 &\int_{\Configs{\Ring{B}}}
   \Poisson{B,\alpha}(A\RConnectedToIn{\xi}C\SymDiff{}C')
   \nu(\xi_{\Ring{B}}=\D{C'})
 \\\stackrel{\text{relax}}{\le}{}
 &\int_{\Configs{\Ring{B}}}
   \Poisson{B,\alpha}(A\RConnectedToIn{\xi}\Complement{B})
   \nu(\xi_{\Ring{B}}=\D{C'})
 \\\stackrel{\eqref{eq_exp_decay_boolean}}{\le}{}
 &K e^{-\kappa{}\SetDistance{A,\Complement{B}}}
 \,.
\end{align*}
For~\eqref{eq_exp_decay_sepev}, assume that $\SetDistance{A,B}>0$.
Choose a sphere $D$ containing $A$ such that \\
$\SetDistance{A,\Ring{D}}>\SetDistance{A,B}$.
Let $D':=\Ring{D}\cup{}B$.
Hence, $\SetDistance{A,D'}\ge\SetDistance{A,B}$, and
\begin{align*}
 &\Modulus{
   \nu(\xi_{A}\in{}E,\xi_{B}\in{}F)
  -\nu(\xi_{A}\in{}E)\nu(\xi_{B}\in{}F)
  }
 \\\stackrel{\eqref{eq_janossy_identities}}{\le}{}
 &\int_{\Configs{D'}}
   \Bigl|
    \HardSphere{D\setminus{}B,C,\lambda}(\xi_{A}\in{}E)
    - \nu(\xi_{A}\in{}E)
   \Bigr|
   \Iverson{C\cap{}B\in{}F}
   \nu(\xi_{D'}=\D{C})
 \\\stackrel{\eqref{eq_exp_decay_fve}}{\le}{}
 &\int_{\Configs{D'}}
   K e^{-\kappa{}\SetDistance{A,D'}}
   \Iverson{C\cap{}B\in{}F}
   \nu(\xi_{D'}=\D{C})
 \\\stackrel{}{=}{}
 &K e^{-\kappa{}\SetDistance{A,B}}
  \nu(\xi_{B}\in{}F)
 \stackrel{}{\le}{}
 K e^{-\kappa{}\SetDistance{A,B}}
 \,.\qedhere
\end{align*}

\end{proof}

\begin{proof}[Proof of Theorem~\ref{thm_exp_decay_rpcf}]
Property~\eqref{eq_dac_twisted_subset} modifies~\eqref{eq_diff_connection_bound} for increasing events to
\begin{equation}
\label{eq_diff_con_bound_inc}
 \Modulus{
  \HardSphere{B,C_1,\lambda}(\xi_{A}\in{}E)
  -\HardSphere{B,C_2,\lambda}(\xi_{A}\in{}E)
 }
 \le
 \Poisson{A,\lambda}(E)\Poisson{B,\lambda}(A\RConnectedToIn{\xi}\Complement{B})
 \,.
\end{equation}
This follows from retracing the second part of the proof of Proposition~\ref{prop_diff_connection_bound} with
\begin{align*}
 \Law(\xi^1_{A}\in{}E,\xi^2_{A}\not\in{}E)
 \stackrel{}{=}{}
 &\Law(\xi^1_{A}\in{}E,\xi^1_{A}\SymDiff\xi^2_{A}\not=\emptyset)
 \\\stackrel{\eqref{eq_dac_def_connected}}{\le}{}
 &\Law(\xi^1_{A}\in{}E,
  A\RConnectedToIn{\xi^1\SymDiff\xi^2}C_1\SymDiff C_2
  )
 \\\stackrel{\text{relax}}{\le}{}
 &\Law(\xi^1_{A}\in{}E,
  A\RConnectedToIn{(\xi^1\cup\xi^2)\setminus{}A}\Complement{B}
  )
 \\\stackrel{\eqref{eq_dac_twisted_subset}}{\le}{}
 &\Law(\xi^3_{A}\in{}E,
  A\RConnectedToIn{\xi^3\setminus{}A}\Complement{B}
  )
 \\\stackrel{\eqref{eq_dac_def_poi}}{=}{}
 &\Poisson{A,\lambda}(E)
  \Poisson{B,\lambda}(A\RConnectedToIn{\xi}\Complement{B})
 \,.
\end{align*}
Retracing the proofs of~\eqref{eq_exp_decay_fve} and~\eqref{eq_exp_decay_sepev} using~\eqref{eq_diff_con_bound_inc} instead of~\eqref{eq_diff_connection_bound} modifies~\eqref{eq_exp_decay_sepev} to:
For all $A,B\in\BoundedBorel$, increasing $E\in\Configs{A}$ and increasing $F\in\Configs{B}$,
\begin{equation}
\label{eq_exp_decay_sepev_inc}
 \Modulus{\nu(\xi_{A}\!\in{}\!E,\xi_{B}\!\in{}\!F)
  -\nu(\xi_{A}\!\in{}\!E)\nu(\xi_{B}\!\in{}\!F)
 }
 \le
  \Poisson{A,\lambda}(E)\Poisson{B,\lambda}(F)
  K e^{-\kappa{}\SetDistance{A,B}}
 \,.
\end{equation}
For all disjoint $A,B\in\BoundedBorel$, bound the second factorial cumulant measure as
\begin{align*}
 &\BigModulus{
  \Expect_\nu\Modulus{\xi_{A}}\Modulus{\xi_{B}}
  - \Expect_\nu\Modulus{\xi_{A}}\Expect_\nu\Modulus{\xi_{B}}
  }
 \\\stackrel{}{\le}{}
 &\sum_{n,m=1}^\infty
  \BigModulus{
    \nu(\Modulus{\xi_{A}}\ge n,\Modulus{\xi_{B}}\ge m)
  - \nu(\Modulus{\xi_{A}}\ge n)\nu(\Modulus{\xi_{B}}\ge m)
  }
 \\\stackrel{\eqref{eq_exp_decay_sepev_inc}}{\le}{}
 &\sum_{n,m=1}^\infty
   \Poisson{A,\lambda}(\Modulus{\xi}\ge n)
   \Poisson{B,\lambda}(\Modulus{\xi}\ge m)
   K e^{-\kappa{}\SetDistance{A,B}}
 \\\stackrel{}{=}{}
 &\lambda\Lebesgue(A)\lambda\Lebesgue(B)K e^{-\kappa{}\SetDistance{A,B}}
 \,.
\end{align*}
Statement~\eqref{eq_exp_decay_rpcf} follows by disintegration with respect to $\lambda^2\Lebesgue^2$.
\end{proof}

\section{Dependently thinning Poisson to hard-sphere}
\label{sec_thin}

\par
Sections~\ref{sec_jj_measure} and~\ref{sec_hs_more} contain additional facts about joint Janossy measures and the hard-sphere model respectively.
Section~\ref{sec_ordering_derivative} describes a measurable total ordering of Euclidean space.
The dependent thinning from a $\Poisson{B,\lambda}$ to a $\HardSphere{B,C,\lambda}$ is in Section~\ref{sec_thin_def}.
This section fixes $\lambda\in[0,\infty[$.
Hence, it drops the $\lambda$ in $\HardSphere{B,C,\lambda}$, $\PartFun{B,C,\lambda}$ and $\Poisson{B,\lambda}$.

\subsection{Joint Janossy measure}
\label{sec_jj_measure}

\par
Let $n\ge 2$ and $\Law$ be a coupling of $n$ PP laws.
A Borel measure $M$ on $(\Configs{B}^n,\ProductAlgebra{B}{n})$ is the \emph{(local) joint Janossy measure} of $\Law$ on $B\in\BoundedBorel$, if, for all $E_1,\dotsc,E_n\in\Algebra{B}$,
\begin{equation}
\label{eq_jj_measure}
 \Law(\forall\,1\le{}i\le{}n: \xi^i_{B}\in{}E_i)
 \\=
  \int_{\Configs{B}^n}
   \prod_{1\le{}i\le{}n}\Iverson{Y_i\in{}E_i}
   M(\D{Y})
 \,.
\end{equation}
Because the local joint Janossy measure on $A\subseteq{}B$ of a coupling $\Law$ on $B$ equals the joint Janossy measure of the restriction of the coupling to $A$, the remainder of this paper drops the quantifier ``local''.
This definition of a joint Janossy measure is between the portmanteau style of the classic case~\eqref{eq_janossy_measure} and the explicit style on generating sets in~\cite[Section 5.3]{Daley_VereJones__AnIntroductionToTheTheoryOfPointProcesses_I__Springer_2003}.
As the sets $\prod_{i=1}^n E_i$ generate $\ProductAlgebra{B}{n}$, there is no loss of generality.
If $\Law$ admits a joint Janossy measure on $B$, then $\Law(\xi_{B}=\D{Y})$ denotes its infinitesimal at $Y\in\Configs{B}^n$.

\par
The identities~\eqref{eq_janossy_identities} generalise directly from the classic to the joint case.
Joint Janossy measures of marginals of a coupling $\Law$ result from integrating out the joint Janossy measure over the complement.

\subsection{More about the hard-sphere model}
\label{sec_hs_more}

\par
The conditional hard-sphere constraint $\SymbolIsHS$ chains.
\begin{equation}
\label{eq_ishs_chains}
 \forall\,X,Y,Z\in\Configs{\RRD}:
 \qquad
 \IsHS{X\cup{}Y|Z} = \IsHS{X|Y\cup{}Z}\IsHS{Y|Z}
 \,.
\end{equation}
The function $\SymbolIsHS$ is $\Configs{\Complement{B}}\times\Configs{B}\to\EnumSet{0,1}$ and measurable on $(\Configs{\Complement{B}}\times\Configs{B},\Algebra{\Complement{B}}\otimes\Algebra{B})$ as a product of measurable functions~\eqref{eq_ishs_def}.
It is monotone decreasing in both arguments.

\par
For $B\in\BoundedBorel$, the function
\begin{equation}
\label{eq_pf_bc_meas_mon}
 \Configs{\Complement{B}}\to[0,\infty[
 \qquad
 C\mapsto\PartFun{B,C}
\end{equation}
is measurable on $(\Configs{\Complement{B}},\Algebra{\Complement{B}})$ and monotone decreasing.
Consequently, $\HardSphere{B,C}$ is measurable in the boundary condition $C$, too.
For $C\in\Configs{\RRD}$, the function
\begin{equation}
\label{eq_pf_dom_mon}
 \BoundedBorel\to[0,\infty[
 \qquad
 B\mapsto\PartFun{B,C\setminus{}B}
\end{equation}
is monotone increasing.
Finally, the relation between~\eqref{eq_hsjan_def} and~\eqref{eq_statmech_partfun} is
\begin{equation}
\label{eq_pf_rewrite}
 \PartFun{B,C}
 =
 \Poisson{B}(\IsHS{\xi|C}=1)\,e^{\lambda\Lebesgue(B)}
 \,.
\end{equation}

\par
The hard-sphere model fulfils the \emph{DLR conditions}~\cite[(2.2)--(2.4)]{Preston__RandomFields__LNM_Springer_1976}.
That is, for $A,B\in\BoundedBorel$, $X\in\Configs{A}$ and $Y\in\Configs{B}$, the Janossy infinitesimal chains.
\begin{equation}
\label{eq_dlr}
 \HardSphere{A\cup{}B,C}(\D{(X\cup{}Y)})
 =
 \HardSphere{A\cup{}B,C}(\xi_{A}=\D{X})
 \HardSphere{B,C\cup{}X}(\D{Y})
 \,.
\end{equation}

\subsection{Ordering and derivative}
\label{sec_ordering_derivative}

\par
This section presents a measurable total ordering of $\RRplusD$.
It allows to define a derivative of measurable functions of Borel subsets of $\RRplusD$.

\par
\begin{figure}
\label{fig_unitsquare}
\begin{center}
\includegraphics[width=0.6\textwidth]{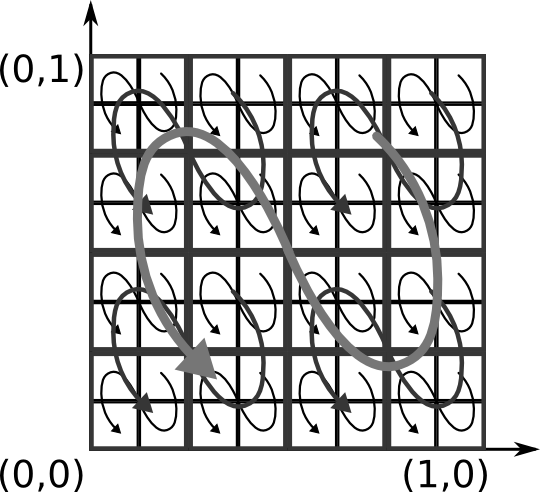}
\end{center}
\caption[Total ordering the unit square]{
Total ordering of the unit square.
The top three levels, with the $i^{th}$ level corresponding the effect of the $i^{th}$ binary digit a point, inducing the ordering are shown.
The first level divides the unit square into four quarter-squares and orders them along the biggest arrow (top-right first, bottom-left last).
Then, within each quarter square, this ordering is repeated (four middle arrows).
The sixteen smallest arrows show the third level, decomposing each sixteenth-square into four parts and ordering them.
The full ordering repeats this recursively on all scales.
}
\end{figure}

The unsigned binary digit sequences are
\begin{equation*}
 \mathcal{D} \
 :=
 \DescSet{\iota\in\EnumSet{0,1}^\ZZ}{\exists{}k:\forall\,n\ge{}k:\iota_n=0}
 \,.
\end{equation*}
A sequence gets assigned a non-negative real value through the map
\begin{equation*}
 b:
 \qquad
 \mathcal{D}\to\RRplus
 \qquad
 \iota\mapsto\sum_{n\in\ZZ}\iota_n 2^n
 \,.
\end{equation*}
There is no inverse of $b$, because multiples of some $2^n$, $n\in\ZZ$ have two preimages under $b$: one ending in an infinite sequence of $0s$ and another one ending in an infinite sequence of $1$s.
Observing that those multiples form a $\Lebesgue$ null-set of $\RRplus$ and choosing the preimage ending in an infinite sequence of $0s$ allows to restrict $b$ to a measurable bijection.
This allows to construct another measurable bijection
\begin{equation*}
 \hat{b}:
 \qquad
 \RRplusD\to\RRplus
 \qquad
 (x_1,\dotsc,x_\Dimension)
 \mapsto
 b^{-1}(
 n\mapsto{}b(x_{(n\operatorname{mod}\Dimension)+1})_{\lfloor{}n/\Dimension\rfloor{}}
 )
 \,.
\end{equation*}
From here on use the bijection $\hat{b}$ implicitly.
See also Figure~\ref{fig_unitsquare}.

\par
The bijection $\hat{b}$ orders $\RRplusD$ measurably and totally.
Denote this order by $\Order$.
The symbols $\pm\infty$ extend $\Order$ with elements being bigger and smaller than each element of $\RRplusD$.
For $a,b\in\RRplusD\cup\EnumSet{\pm\infty}$ with $a\LT{}b$, there is the \emph{interval} $]a,b]:=\DescSet{x}{a\LT x\LTE b}$, as well as all standard variations thereof.

\par
Let $B\in\BoundedBorel$ with $B\subseteq\RRplusD$ and $\Lebesgue(B)>0$.
As $\Lebesgue$--a.e. $x\in{}B$ is a density point of $\Lebesgue$~\cite[Section 5.8(ii)]{Bogachev__MeasureTheory__Springer_2007}, there exists $\varepsilon_x>0$, such that for all $0<\varepsilon<\varepsilon_x$ there exist points $x_\varepsilon^-,x_\varepsilon^+\in{}B$ with  $x_\varepsilon^-\LT{}x\LT{}x_\varepsilon^+$ such that $\Lebesgue(]x_\varepsilon^-,x[)=\Lebesgue(]x,x_\varepsilon^+[)=\varepsilon$.
The derivative of $f:{}]a,b[\to\RR$ at $x$ is
\begin{equation}
\label{eq_def_derivative}
 f'(x)
 :=
 \lim_{\varepsilon\to 0}
  \frac{f(x_\varepsilon^+)-f(x_\varepsilon^-)}{2\varepsilon}
 \,,
\end{equation}
whenever this limit is defined.
This is the usual one-dimensional derivative on $\RRplus$ mapped back through $\hat{b}^{-1}$.

\subsection{The thinning}
\label{sec_thin_def}

\par
This section presents a coupling between a hard-sphere PP and a dominating Poisson PP.
The coupling is an explicit dependent thinning from the dominating Poisson PP.
The thinning probability\footnote{The established name for keeping a point as part of the smaller thinned PP.} is related to the logarithm of the free energy.
Its explicit form implies the measurability of the coupling with respect to the boundary condition.

\par
For the remainder of this section, fix $B\in\BoundedBorel$ and $C\in\Configs{\Complement{B}}$.
Without loss of generality, translation-invariance of the hard-sphere model lets us restrict to $B\subseteq\RRplusD$.
This way, the order from Section~\ref{sec_ordering_derivative} applies.
For the remainder of Section~\ref{sec_thin}, restrict intervals to $B$, i.e., $]a,b]$ denotes $]a,b]\cap{}B$.

\begin{proposition}
\label{prop_derivative}
For $\Lebesgue$--a.e. $x\in{}B$ and $Y\in\Configs{]-\infty,x[}$,
\begin{equation}
\label{eq_derivative}
 -\frac{1}{\lambda}
  \frac{\partial}{\partial x}\log\PartFun{]x,\infty[,C\cup{}Y}
 \stackrel{}{=}
 \IsHS{\EnumSet{x}|C\cup{}Y}
 \frac%
  {\PartFun{]x,\infty[,C\cup{}Y\cup\EnumSet{x} }}
  {\PartFun{]x,\infty[,C\cup{}Y}}
 \,.
\end{equation}
\end{proposition}

\par
The proof of Proposition~\ref{prop_derivative} is in Section~\ref{sec_prop_derivative_proof}.
Proposition~\ref{prop_derivative} calculates the derivative of the free energy of a right-unbounded interval in $B$ with respect to the Poisson intensity.
The monotonicity of $\SymbolPartFun$ in the domain~\eqref{eq_pf_dom_mon} applied to the lhs of~\eqref{eq_derivative} implies its monotone growth in $x$ outside of $\Sphere{C\cup{}Y}$.
The monotonicity applied to the rhs of~\eqref{eq_derivative} implies that its value lies in $[0,1]$.
For each $Y\in\Configs{B}$, this yields a restricted thinning kernel
$\OneShot{Y}:{}]\max{}Y,\infty[\to[0,1]$, with $\OneShot{Y}(x)$ given by the lhs of~\eqref{eq_derivative} and interpreted as the probability of keeping the point $x$.

\par
The thinning arises from an ordered exploration of $B$.
It explores the points of a Poisson realisation $Y_2\in\Configs{B}$ in the order induced by $\Order$.
The starting thinning kernel is $\OneShot{\emptyset}$.
At the first point $y_1\in{}Y_2$ it keeps, it is replaced by the kernel  $\OneShot{\EnumSet{y_1}}$.
The second kernel explores $Y_2\cap{}]y_1,\infty[$.
At the first point $y_2\in{}Y_2\cap{}]y_1,\infty[$ it keeps, it is replaced by the kernel $\OneShot{\EnumSet{y_1,y_2}}$ exploring $Y_2\cap{}]y_2,\infty[$.
Iterate until all of $Y_2$ has been explored.
The usage of the restricted thinning kernels is ``one-shot'', i.e., a kernel is used until the first time it keeps a point from $Y_2$ and then replaced by the next thinning kernel with updated dependencies.
The following definition formalises this dependent update of the thinning kernel.

\begin{definition}
\label{def_thin}
\begin{subequations}
\label{eq_thin_def}
For $x\in{}B$ and $Y_1\in\Configs{B}$, the dependent thinning probability is
\begin{equation}
\label{eq_thin_proba}
 \ThinProba{x}{Y_1} := \OneShot{Y_1\cap{}]-\infty,x[}(x)\,.
\end{equation}
The choice function distinguishes between kept and deleted points.
\begin{equation}
\label{eq_thin_choose}
 \ThinChoice{x,Y_1} :=
 \Iverson{x\in{}Y_1}\ThinProba{x}{Y_1}
 + \Iverson{x\not\in{}Y_1}(1-\ThinProba{x}{Y_1})
 \,.
\end{equation}
The Janossy infinitesimal of the thinning is
\begin{equation}
\label{eq_thin_janossy}
 \HSPoiThin{B,C}(\D{Y})
 :=
 \Iverson{Y_1\subseteq{}Y_2}
 \left(
  \prod_{x\in{}Y_2} \ThinChoice{x,Y_1}
 \right)
 \Poisson{B}(\D{Y_2})
 \,.
\end{equation}
\end{subequations}
\end{definition}

The dependent thinning kernel $\ThinProba{.}{Y_1}$ is a piece-wise combination of the restricted kernels.
The term $Y_1\cap{}]-\infty,x[$ in the rhs of~\eqref{eq_thin_proba} selects the appropriate restricted kernel, depending on the already explored and kept points.
The choice function~\eqref{eq_thin_choose} assigns correct probabilities to points being kept as part of $Y_1$ and deleted on $B\setminus{}Y_1$ respectively.
Finally, the Janossy infinitesimal~\eqref{eq_thin_janossy} describes the joint probability of a Poisson realisation $Y_2$ and keeping exactly the subset $Y_1$ for the thinned process.

\begin{theorem}
\label{thm_thin}
The dependent thinning $\HSPoiThin{B,C}$ is a dominating coupling between a Poisson PP and a hard-sphere PP, as
\begin{subequations}
\label{eq_thin_properties}
\begin{gather}
 \label{eq_thin_hs}
 \HSPoiThin{B,C}(\xi^1=\D{Y}) = \HardSphere{B,C}(\D{Y})
 \,,
 \\
 \label{eq_thin_poi}
 \HSPoiThin{B,C}(\xi^2=\D{Y}) = \Poisson{B}(\D{Y})
\end{gather}
and
\begin{equation}
\label{eq_thin_subset}
 \HSPoiThin{B,C}(\xi^1\subseteq{}\xi^2) = 1
 \,.
\end{equation}
\end{subequations}
The boundary condition may be restricted to $\Ring{B}$.
The law $\HSPoiThin{B,C}$ is measurable in $C$.
\end{theorem}

The proof of Theorem~\ref{thm_thin} is in Section~\ref{sec_thm_thin_proof}.

\subsection{Proof of Proposition~\ref{prop_derivative}}
\label{sec_prop_derivative_proof}

\par
As $\Lebesgue$--a.e. $\Lebesgue(]\max Y,x])>0$, rephrase the dependence between $x$ and $Y$ in~\eqref{eq_derivative}.
For $\Lebesgue$--a.e. $a\in{}B\cup\EnumSet{-\infty}$, $b\in{}]a,\infty[$, $Y\in\Configs{]-\infty,a]}$ and $x\in{}]a,b[$, the aim is to show that~\eqref{eq_derivative} holds.
Let $C':=C\cup{}Y$.
Regard the measurable functions
\begin{equation}
\label{eq_interval_func_1}
\begin{aligned}
 h:\qquad
 &]a,b]\to[0,1]
 &x
 &\mapsto\IsHS{\EnumSet{x}|C'}
 \,,\\
 z:\qquad
 &[a,b]\to[0,\infty[
 &x
 &\mapsto\PartFun{]x,\infty[,C'}
 \,,\\
 s:\qquad
 &[a,b]\to[0,\infty[
 &x
 &\mapsto\PartFun{]x,\infty[,C'\cup\EnumSet{x}}
 \,.
\end{aligned}
\end{equation}
Using the derivative~\eqref{eq_def_derivative}, if $z'=-\lambda h s$ on $]a,b[$, then~\eqref{eq_derivative} follows from
\begin{equation*}
 -\left(\frac{\log z}{\lambda}\right)'
 = -\frac{z'}{\lambda z}
 = -\frac{-\lambda h s}{\lambda z}
 = \frac{h s}{z}
 \,.
\end{equation*}
The remainder of this section shows that $z'=-\lambda h s$ $\Lebesgue$--a.e. on $]a,b[$.

\par
Wlog assume that $\Lebesgue(]a,b[)>0$.
Fix $x\in{}]a,b[$.
Using the notation from~\eqref{eq_def_derivative}, let $A_\varepsilon:={}]x_\varepsilon^-,x_\varepsilon^+]$, for $\varepsilon<\varepsilon_x$.
If $\varepsilon$ is small enough, then $A_\varepsilon\subseteq\Sphere{y}$ holds uniformly in  $y\in{}A_\varepsilon$.
Let $A_\varepsilon^-:={}]x_\varepsilon^-,\infty[$ and $A_\varepsilon^+:=[x_\varepsilon^+,\infty[{}=A_\varepsilon^-\setminus{}A_\varepsilon$.
Using~\eqref{eq_pf_rewrite}, expand $z(x_\varepsilon^-)$ and $z(x_\varepsilon^+)$ as
\begin{align*}
 z(x_\varepsilon^-)
 \stackrel{}{=}{}
 &e^{\lambda\Lebesgue(A_\varepsilon^-)}
 \int_{\Configs{A_\varepsilon^-}}
 \IsHS{Z|C'}
 \Poisson{A_\varepsilon^-}(\D{Z})
 \\\stackrel{}{=}{}
 &e^{\lambda\Lebesgue(A_\varepsilon^-)}
 \int_{\Configs{A_\varepsilon}}
 \IsHS{X|C'}
 \int_{\Configs{A_\varepsilon^+}}
 \IsHS{Z|C'\cup{}X}
 \Poisson{A_\varepsilon^+}(\D{Z})
 \Poisson{A_\varepsilon}(\D{X})
 \,,\\
 z(x_\varepsilon^+)
 \stackrel{}{=}{}
 &e^{\lambda\Lebesgue(A_\varepsilon^+)}
 \int_{\Configs{A_\varepsilon^+}}
 \IsHS{Z|C'}
 \Poisson{A_\varepsilon^+}(\D{Z})
 \\\stackrel{}{=}{}
 &e^{\lambda\Lebesgue(A_\varepsilon^-)}
 \int_{\Configs{A_\varepsilon}}
 \Iverson{X=\emptyset}
 \IsHS{X|C'}
 \int_{\Configs{A_\varepsilon^+}}
 \IsHS{Z|C'\cup{}X}
 \Poisson{A_\varepsilon^+}(\D{Z})
 \Poisson{A_\varepsilon}(\D{X})
 \,,
\end{align*}
to see that $z'(x)=\displaystyle\lim_{\varepsilon\to 0}
 \frac{z(x_\varepsilon^+)-z(x_\varepsilon^-)}{2\varepsilon}$ equals
\begin{equation*}
 \lim_{\varepsilon\to 0}
  \frac{e^{\lambda\Lebesgue(A_\varepsilon^-)}}{2\varepsilon}
  \int_{\Configs{A_\varepsilon}}
  -\Iverson{X\not=\emptyset}
  \IsHS{X|C'}
  \int_{\Configs{A_\varepsilon^+}}
  \IsHS{Z|C'\cup{}X}
  \Poisson{A_\varepsilon^+}(\D{Z})
  \Poisson{A_\varepsilon}(\D{X})
 \,.
\end{equation*}
The case $\Cardinality{X}\ge{}2$ is irrelevant, because all integrands take values in $[-1,1]$ and $\Poisson{A_\varepsilon}(\Cardinality{\xi}\ge{}2) = o(\varepsilon^2)$.
Hence, in the case $\Cardinality{X}=1$, let $y$ be the single point in $X$ and rewrite $z'(x)$ into
\begin{equation*}
 -\lim_{\varepsilon\to 0}
  \frac{e^{\lambda\Lebesgue(A_\varepsilon^-)}}{2\varepsilon}
  \int_{A_\varepsilon}
  \IsHS{\EnumSet{y}|C'}
  \int_{\Configs{A_\varepsilon^+}}
  \IsHS{Z|C'\cup\EnumSet{y}}
  \Poisson{A_\varepsilon^+}(\D{Z})
  e^{-2\lambda\varepsilon}
  \lambda\D{y}
 \,.
\end{equation*}
Expand the domain of the inner integration from $A_\varepsilon^+$ to $]y,\infty[$ paying a penalty of $e^{\lambda\Lebesgue(]y,x_\varepsilon^+])}$.
As $\Lebesgue(A_\varepsilon^-) + \Lebesgue(]y,x_\varepsilon^+]) - 2\varepsilon = \Lebesgue(]y,\infty[)$, rewrite $z'(x)$ into
\begin{equation*}
 -\lim_{\varepsilon\to 0}
  \frac{1}{2\varepsilon}
  \int_{A_\varepsilon}
  \underbrace{
  \IsHS{\EnumSet{y}|C'}
  }_{=h(y)}
  \underbrace{
  e^{\lambda\Lebesgue(]y,\infty[)}
  \int_{\Configs{]y,\infty[}}
  \IsHS{Z|C'\cup\EnumSet{y}}
  \Poisson{]y,\infty[}(\D{Z})
  }_{=s(y)\text{ by~\eqref{eq_pf_rewrite}}}
  \lambda\D{y}
 \,.
\end{equation*}
The \emph{Lebesgue differentiation theorem}~\cite[Thm 5.6.2]{Bogachev__MeasureTheory__Springer_2007} implies that, $\Lebesgue$--a.e.,
\begin{equation*}
 z'(x)
 \stackrel{}{=}{}
 -\lim_{\varepsilon\to 0}
  \frac{1}{2\varepsilon}
  \int_{A_\varepsilon}
  h(y)
  s(y)
  \lambda\D{y}
 =
 -\lambda h(x)s(x)
 \,.
\end{equation*}

\subsection{Deleting all points and an integral equation}

This section calculates the probability of deleting all points within an interval of $B$.
The key point is that the dependent thinning kernel reduces to a single restricted thinning kernel.
For all $a,b\in{}B\cup\EnumSet{\pm\infty}$ with $a\LT b$, $Y\in\Configs{]-\infty,a]}$ and $X\in\Configs{[b,\infty[}$,
\begin{equation}
\label{eq_thin_interval}
 \int_{\Configs{]a,b[}}
  \prod_{z\in{}Z} \ThinChoice{z,Y\cup{}X}
  \Poisson{]a,b[}(\D{Z})
 =
 \frac{\PartFun{]b,\infty[,C\cup{}Y}}
      {\PartFun{]a,\infty[,C\cup{}Y}}
 \,.
\end{equation}
The case $]-\infty,\infty[{}=B$ implies $Y=\emptyset$ and the correct probability $1/\PartFun{B,C}$.
\par
The solution of~\eqref{eq_thin_interval} comes from an integral equation.
For each $x\in{}]a,b[\,$, as $(Y\cup{}X)\cap{}]-\infty,x[{}=Y$, $\ThinChoice{z,Y\cup{}X}=1-\ThinProba{x}{Y}=1-\OneShot{Y}(x)$.
Regard the measurable functions
\begin{equation}
\label{eq_interval_func_2}
\begin{aligned}
 q:\qquad
 &]a,b[\to[0,1]
 &x
 &\mapsto{}1-\OneShot{Y}(x)
 \,,\\
 l:\qquad
 &[a,b]\to[0,\infty[
 &x
 &\mapsto\int_{]x,b[}1dy = \Lebesgue(]x,b[)
 \,,\\
 e:\qquad
 &[a,b]\to[1,\infty[
 &x
 &\mapsto{}e^{\lambda\Lebesgue(]x,b[)}=e^{\lambda l(x)}
 \,,\\
 t:\qquad
 &[a,b]\to[0,1]
 &x
 &\mapsto
  e(x)
  \int_{\Configs{]x,b[}}
   \prod_{z\in{}Z} q(z)
   \Poisson{]x,b[}(\D{Z})
 \,.
\end{aligned}
\end{equation}
Showing~\eqref{eq_thin_interval} is equivalent to calculating $t(a)/e(a)$.
If $]x,b[$ contains a point, then splitting the smallest point off yields an integral equation for $t$.
\begin{subequations}
\label{eq_inteq}
\begin{equation}
\label{eq_inteq_eq}
\begin{aligned}
 t(x)
 \stackrel{}{=}{}
 &e(x)\int_{\Configs{]x,b[}}
   \prod_{z\in{}Z} q(z)
   \Poisson{]x,b[}(\D{Z})
 \\\stackrel{}{=}{}
 &e(x)\Poisson{]x,b[}(\xi=\emptyset)
 \\&+
  e(x)\int_{]x,b[}
   q(y)
   e^{-\lambda\Lebesgue(]x,y[)}
   \int_{\Configs{]y,b[}}
    \prod_{z\in{}Z} q(z)
    \Poisson{]y,b[}(\D{Z})
  \lambda{}\D{y}
 \\\stackrel{}{=}{}
 &1 + \lambda \int_{]x,b[} q(y) t(y) \D{y}
 \,.
\end{aligned}
\end{equation}
There is a boundary condition
\begin{equation}
\label{eq_inteq_bc}
 t(b)
 = e(b)\int_{\Configs{\emptyset}} 1 \Poisson{\emptyset}(\D{Z})
 = 1
 \,.
\end{equation}
\end{subequations}
Thus, a solution of~\eqref{eq_inteq} yields~\eqref{eq_thin_interval}.
Because the setup in~\eqref{eq_interval_func_1} is the same as in~\eqref{eq_interval_func_2}, consider
\begin{align*}
 g:
 \qquad
 [a,b]\to[0,\infty[
 \qquad
 x\mapsto\frac{e(x)}{z(x)}
 \,.
\end{align*}
Establish that $l, e, z$ en $g$ are absolutely continuous by showing that they are Lipschitz continuous on $[a,b]$~\cite[Lemma 5.3.2]{Bogachev__MeasureTheory__Springer_2007}.
Both $e$ and $z$ are monotone decreasing with bounds $e(a)\ge{}e(x)$ and $z(a)\ge{}z(x)\ge{}z(b)\ge 1$, for $x\in[a,b]$.
Because $l$ has Lipschitz constant $1$, $e$ has Lipschitz constant $\lambda{}e(a)$.
For $a\LTE{}x\LT{}y\LTE{}b$, use~\eqref{eq_statmech_partfun} to expand $z$ and obtain the bound
\begin{multline*}
 \Modulus{z(x)-z(y)}
 \le
 \sum_{n=1}^\infty \frac{\lambda^n}{n!}
  \BigModulus{
    \int_{]x,\infty[^n\setminus{}]y,\infty[^n}
      e^{-H(a_1,\dotsc,a_n|C)} \prod_{i=1}^n \D{a_i}
    }
 \\
 \le
 \sum_{n=1}^\infty \frac{\lambda^n}{n!}
  \Lebesgue(]x,\infty[^n\setminus{}]y,\infty[^n)
 \le
 \sum_{n=1}^\infty \frac{\lambda^n}{n!} n \Modulus{x-y} \Lebesgue(]y,\infty[)^{n-1}
 =
 \lambda e(y) \Modulus{x-y}
 \,.
\end{multline*}
Hence, $z$ has Lipschitz constant $\lambda e(a)$.
For $a\LTE{}x\LT{}y\LTE{}b$, bound $g$ by
\begin{multline*}
 \Modulus{g(x)-g(y)}
 =
 \BigModulus{\frac{e(x)z(y)-e(y)z(y)+e(y)z(y)-e(y)z(x)}{z(x)z(y)}}
 \\
 \le
 \frac{\Modulus{e(x)-e(y)}}{z(x)} + \frac{e(y)\Modulus{z(x)-z(y)}}{z(x)z(y)}
 \,.
\end{multline*}
Applying the bounds for $e$ and $z$, $g$ has Lipschitz constant $\lambda\frac{e(a)}{z(b)}+\lambda\frac{e(a)^2}{z(b)^2}$.

\par
As the functions $l, e, z$ en $g$ are all absolutely continuous with respect to $\Lebesgue$, they are Lebesgue differentiable $\Lebesgue$--a.e.~.
Because $l'=-1$ and $e'=-\lambda e$, $\Lebesgue$--a.e. on $]a,b[$,
\begin{equation*}
 g'
 = \frac{ze'-ez'}{z^2}
 =
 \frac{z(-\lambda e)- e(-\lambda hs)}{z^2}
 = -\frac{\lambda{}e}{z}\left(1-\frac{hs}{z}\right)
 \stackrel{\eqref{eq_derivative}}{=}
 -\lambda g q
 \,.
\end{equation*}
As $e(b)=1$, $g(b)=\frac{1}{z(b)}\le 1$, by~\eqref{eq_pf_bc_meas_mon}.
Integration of $g$ yields
\begin{equation*}
 g(x)
 = g(b) - \int_{]x,b[} g'(y) \D{y}
 = \frac{1}{z(b)} + \lambda \int_{]x,b[} g(y) q(y) \D{y}
 \,.
\end{equation*}
Hence, the function $z(b)g$ solves~\eqref{eq_inteq} $\Lebesgue$--a.e. and is $\Lebesgue$--a.e. equal to $t$.
The expression $\frac{t(a)}{e(a)}=\frac{z(b)g(a)}{e(a)}=\frac{z(b)}{e(a)}\frac{e(a)}{z(a)}=\frac{z(b)}{z(a)}$ yields the rhs of~\eqref{eq_thin_interval}.

\subsection{Proof of Theorem~\ref{thm_thin}}
\label{sec_thm_thin_proof}

\par
The thinning $\HSPoiThin{B,C}$ is well-defined, as the choice function sums to one over all possible choices.
\begin{equation}
\label{eq_thin_sum_prod}
 \forall\,\text{ disjoint }X,Z\in\Configs{B}:
 \qquad
 \beta(X,Z) :=
 \sum_{Y\subseteq{}Z} \prod_{x\in{}Z} \ThinChoice{x,X\cup{}Y} = 1
 \,.
\end{equation}
Ascertain~\eqref{eq_thin_sum_prod} by induction on the size of $Z$.
The base case $Z=\emptyset$ is trivially true.
Otherwise, let $z:=\min Z$ (with respect to $\Order$) and $Z':=Z\setminus\EnumSet{z}$.
For each $Y\subseteq{}Z$, the fact that $Y\subseteq[z,\infty[$ implies that $(X\cup{}Y)\cap]-\infty,z]=X\cap]-\infty,z]$, whence
\begin{equation}
\label{eq_thin_sum_prod_simplify}
   \ThinChoice{z,X\cup{}Y}
 = \ThinProba{z}{X\cup{}Y}
 = \OneShot{X\cap{}]-\infty,z[}(z)
 \,.
\end{equation}
Hence, with $P:= \OneShot{X\cap{}]-\infty,z[}(z)$,
\begin{align*}
 \beta(X,Z)\stackrel{\eqref{eq_thin_sum_prod}}{=}{}
 &\sum_{Y'\subseteq{}Z'}
  \left(
  \prod_{x\in{}Z} \ThinChoice{x,X\cup\EnumSet{z}\cup{}Y'}
  + \prod_{x\in{}Z} \ThinChoice{x,X\cup{}Y'}
  \right)
 \\\stackrel{\eqref{eq_thin_sum_prod_simplify}}{=}{}
 &\sum_{Y'\subseteq{}Z'}
   P
   \prod_{x\in{}Z'} \ThinChoice{x,X\cup\EnumSet{z}\cup{}Y'}
 +\sum_{Y'\subseteq{}Z'}
   (1-P)\prod_{x\in{}Z'} \ThinChoice{x,X\cup{}Y'}
 \\\stackrel{\eqref{eq_thin_sum_prod}}{=}{}
 &P\beta(X\cup{}\EnumSet{z},Z') + (1-P) \beta(X,Z')
 \stackrel{\text{ind}}{=}{}
 P + 1 - P
 = 1
 \,.
\end{align*}

\par
The freedom to restrict the boundary condition $C$ to $C\cap\Ring{B}\in\Configs{\Ring{B}}$ follows from the same freedom for $\SymbolIsHS$ and $\SymbolPartFun$~\eqref{eq_statmech_partfun}.
The measurability in the boundary condition follows from the measurability preserving operations in~\eqref{eq_thin_def} and the measurability of the rhs of~\eqref{eq_derivative}.

\par
The fact that $\HSPoiThin{B,C}$ is a stochastic domination is evident from the construction.
The construction as thinning implies that the second marginal is Poisson~\eqref{eq_thin_poi}.
The remainder of this section shows that the first marginal is hard-sphere~\eqref{eq_thin_hs}.

\par
Let $Y\in\Configs{B}$ with $n:=\Cardinality{Y}$.
Order $Y=:\EnumSet{y_1,\dotsc,y_n}$ increasingly by $\Order$.
Let $y_0:=-\infty$ and $y_{n+1}:=\infty$.
For $0\le{}i\le{}n$, let $B_i:={}]y_i,y_{i+1}[$, $A_i:={}]y_i,\infty[$, $Y_i:=\EnumSet{y_1,\dotsc,y_i}$ and $C_i:=C\cup{}Y_i$.
Thus, $Y\cap{}]-\infty,y_i[{}=Y_{i-1}$.
\begin{align*}
 &\HSPoiThin{B,C}(\xi^1=\D{Y})
 \\\stackrel{\eqref{eq_thin_janossy}}{=}{}
 &\int_{\Configs{B}}
   \Iverson{Y\subseteq{}Z}
   \prod_{z\in{}Z} \ThinChoice{z,Y}
   \Poisson{B}(\D{Z})
 \\={}
 &\left(
  \int_{\Configs{B\setminus{}Y}}
   \prod_{z\in{}Z} \ThinChoice{z,Y}
   \Poisson{B}(\D{Z})
  \right)
  \left(\prod_{y\in{}Y} \ThinChoice{y,Y}\right)
  e^{\lambda\Lebesgue(B)}
  \Poisson{B}(\D{Y})
 \\\stackrel{\eqref{eq_thin_choose}}{=}{}
 &\left(
  \int_{\Configs{B\setminus{}Y}}
   \prod_{z\in{}Z} \ThinChoice{z,Y}
   \Poisson{B\setminus{}Y}(\D{Z})
  \right)
  \left(\prod_{y\in{}Y} \ThinProba{y}{Y}\right)
  e^{\lambda\Lebesgue(B)}
  \Poisson{B}(\D{Y})
 \\\stackrel{\eqref{eq_thin_proba}}{=}{}
 &\left(
   \prod_{i=0}^n
   \int_{\Configs{B_i}}
    \prod_{z\in{}Z} \ThinChoice{z,Y}
    \Poisson{B_i}(\D{Z})
  \right)
  \left(
  \prod_{i=1}^n
   \OneShot{Y_{i-1}}(y_i)
  \right)
  e^{\lambda\Lebesgue(B)}
  \Poisson{B}(\D{Y})
 \,.
\end{align*}
For $0\le{}i\le{}n$,
\begin{equation*}
 \int_{\Configs{B_i}}
  \prod_{z\in{}Z} \ThinChoice{z,Y}
  \Poisson{B_i}(\D{Z})
 \stackrel{\eqref{eq_thin_interval}}{=}{}
 \frac{\PartFun{A_{i+1},C_{i}}}
      {\PartFun{A_{i},C_{i}}}
 \,.
\end{equation*}
For $1\le{}i\le{}n$,
\begin{align*}
 \OneShot{Y_{i-1}}(y_i)
 \stackrel{\eqref{eq_derivative}}{=}{}
 \IsHS{\EnumSet{y_i}|C_{i-1}}
 \frac{\PartFun{A_{i},C_{i}}}
      {\PartFun{A_{i},C_{i-1}}}
 \,.
\end{align*}
Combine these rewritings to see that $\HSPoiThin{B,C}(\xi^1=\D{Y})$ equals
\begin{equation*}
  \left(
   \prod_{i=0}^n
   \frac{\PartFun{A_{i+1},C_{i}}}{\PartFun{A_{i},C_{i}}}
  \right)
  \left(
   \prod_{i=1}^n
   \IsHS{\EnumSet{y_i}|C_{i-1}}
   \frac{\PartFun{A_{i},C_{i}}}{\PartFun{A_{i},C_{i-1}}}
  \right)
  e^{\lambda\Lebesgue(B)}
  \Poisson{B}(\D{Y})
  \,.
\end{equation*}
Combine the hard-sphere constraints by~\eqref{eq_ishs_chains}.
Join the two products and cancel the factors except the denominator at index $0$ and the numerator at index $n$ from the left product.
\begin{multline*}
 \HSPoiThin{B,C}(\xi^1=\D{Y})
 ={}
 \frac{\PartFun{A_{n+1},C_n}}
       {\PartFun{A_0,C_0)}}
  \IsHS{Y|C}
  e^{\lambda\Lebesgue(B)}
  \Poisson{B}(\D{Y})
 \\\stackrel{}{=}{}
 \frac{\PartFun{\emptyset,C\cup{}Y}\IsHS{Y|C}\Poisson{B}(\D{Y})}
       {\PartFun{B,C}e^{-\lambda\Lebesgue(B)}}
 \stackrel{\eqref{eq_pf_rewrite}}{=}{}
 \frac{\IsHS{Y|C}\Poisson{B}(\D{Y})}
       {\Poisson{B}(\D{Y}|\IsHS{\xi|C}=1)}
 \stackrel{\eqref{eq_hsjan_def}}{=}{}
 \HardSphere{B,C}(\D{Y})
 \,.
\end{multline*}

\section{The twisted coupling family}
\label{sec_twisted}

Definition~\ref{def_twisted} defines a family of couplings recursively.
Proposition~\ref{prop_twisted_rec} shows that it is a disagreement coupling family of intensity $\lambda$ for the hard-sphere model.
The notational conventions outlined at the beginning of Section~\ref{sec_thin} apply.

\begin{definition}
\label{def_twisted}
\begin{subequations}
\label{eq_twisted_def}
Let $B\in\BoundedBorel$ and $C_1,C_2\in\Configs{\Complement{B}}$.
For $1\le{}i\le{}2$, let $F_i := B\cap\Ring{C_i}$.
Let $D:=F_1\cup{}F_2$ be the \emph{zone of disagreement} and partition it into $D_1:=F_2\setminus{}F_1$, $D_2:=F_1\setminus{}F_2$ and $D_0:=F_1\cap{}F_2$.\\
Define the joint Janossy intensity of the law $\TwistedZone{B,C_1,C_2}$ on $(\Configs{D}^3,\ProductAlgebra{D}{3})$ by
\begin{equation}
\label{eq_twisted_def_zone}
\begin{aligned}
 \TwistedZone{B,C_1,C_2}(\D{Y})
 :={}
 &\Poisson{D_0}(\D{(Y_3\cap{}D_0)})
 \\\times{}
 &\Iverson{Y_1\subseteq{}D_1}
 \HSPoiThin{B,C_1}(\xi_{D_1}=\D{(Y_1\cap{}D_1,Y_3\cap{}D_1)})
 \\\times{}
 &\Iverson{Y_2\subseteq{}D_2}
 \HSPoiThin{B,C_2}(\xi_{D_2}=\D{(Y_2\cap{}D_2,Y_3\cap{}D_2)})
 \,.
\end{aligned}
\end{equation}
Define the joint Janossy intensity of law $\TwistedRec{B,C_1,C_2}$ on $(\Configs{B}^3,\ProductAlgebra{B}{3})$ recursively by
\begin{multline}
\label{eq_twisted_def_rec}
 \TwistedRec{B,C_1,C_2}(\D{Y})
 :=
 \Iverson{D=\emptyset}
 \Iverson{Y_1=Y_2}
 \HSPoiThin{B,C_1\cup{}C_2}(\D{(Y_1,Y_3)})
 \\+
 \Iverson{D\not=\emptyset}
 \TwistedZone{B,C_1,C_2}(\D{(Y\cap{}D)})
 \TwistedRec{B\setminus{}D,Y_1\cap{}D,Y_2\cap{}D}(\D{(Y\setminus{}D)})
 \,.
\end{multline}
\end{subequations}
\end{definition}

\par
The idea behind the recursive construction of $\TwistedRec{}$ is as follows:
The sets $F_1$ and $F_2$ describe the parts of the domain forbidden by the respective boundary conditions.
If one can construct the disagreement coupling on $D$, then recursion takes care of the rest.
\par
If $D=\emptyset$, a dominating coupling with an identification of the two hard-sphere PPs is already a disagreement coupling.
Since there are no disagreeing points, no connection to the disagreeing boundary is needed.
For $1\le{}i\le{}2$, let $C_i':=C_i\cap\Ring{B}$.
Hence, $C_1'\SymDiff{}C_2'=\emptyset$ and $C_1'\cup{}C_2'=C_1'=C_2'$.
\par
If $D\not=\emptyset$, then the partition $\EnumSet{D_0,D_1,D_2}$ comes into play.
Points of $\xi^1$ and $\xi^2$ can only lie in $D_1$ and $D_2$ respectively.
Independent projections of a dominating coupling take care of that.
This also connects the disagreeing points to the boundary for free.
On $D_0$, an independent Poisson PP of intensity $\lambda$ ensures that there is a Poisson PP on all of $D$.
This is the ``twist''.

\par
The event of connecting disagreement with the boundary in $\ProductAlgebra{B}{2}$ is
\begin{equation}
\label{eq_def_disagreement_connection}
 \DisagreementConnection{B}{C_1,C_2}
 :=
 \DescSet{Y\in\Configs{B}^2}{
  \forall\,x\in{}Y_1\SymDiff{}Y_2:
  x\RConnectedToIn{Y_1\SymDiff{}Y_2}C_1\SymDiff{}C_2
 }
 \,.
\end{equation}

\begin{proposition}
\label{prop_twisted_zone}
\begin{subequations}
The boundary conditions of $\TwistedZone{B,C_1,C_2}$ may be restricted to $\Configs{\Ring{B}}$.
The law $\TwistedZone{B,C_1,C_2}$ is jointly measurable in $(C_1,C_2)$ and has the right marginals
\begin{equation}
\label{eq_twisted_zone_hs}
 \forall\,1\le{}i\le{}2:\qquad
 \TwistedZone{B,C_1,C_2}(\xi^i=\D{Y})
 = \HardSphere{B,C_i}(\xi_{D}=\D{Y})
 \,,
\end{equation}
\begin{equation}
\label{eq_twisted_zone_poi}
 \TwistedZone{B,C_1,C_2}(\xi^3=\D{Y})
 = \Poisson{D}(\D{Y})
 \,.
\end{equation}
It also has the useful properties
\begin{equation}
\label{eq_twisted_zone_contained}
 \forall\,1\le{}i\le{}2:\qquad
 \TwistedZone{B,C_1,C_2}(\xi^i\subseteq{}D_i) = 1
 \,,
\end{equation}
\begin{equation}
\label{eq_twisted_zone_intersection_empty}
 \TwistedZone{B,C_1,C_2}(\xi^1\cap\xi^2=\emptyset) = 1
 \,,
\end{equation}
\begin{equation}
\label{eq_twisted_zone_dominate}
 \forall\,1\le{}i\le{}2:\qquad
 \TwistedZone{B,C_1,C_2}(\xi^i\subseteq\xi^3) = 1
 \,,
\end{equation}
\begin{equation}
\label{eq_twisted_zone_connect}
 \TwistedZone{B,C_1,C_2}(
  (\xi^1,\xi^2)\in\DisagreementConnection{D}{C_1,C_2}
 ) = 1
 \,.
\end{equation}
\end{subequations}
\end{proposition}

\begin{proof}
\par
The freedom to restrict the boundary conditions to $\Configs{\Ring{B}}$ follows from the same property of $\HSPoiThin{B,.}$ in Theorem~\ref{thm_thin}.
The measurability in the boundary conditions follow from the same properties of the law $\Poisson{B}$ and $\HSPoiThin{B,.}$ in Theorem~\ref{thm_thin} and the other measurable indicator terms in construction~\eqref{eq_twisted_def_zone}.

\par
The Poisson marginal~\eqref{eq_twisted_zone_poi} is a straightforward integration over~\eqref{eq_twisted_def_zone}.
The hard-sphere marginals~\eqref{eq_twisted_zone_hs} use the hard-core exclusion together with the properties of the partition $\EnumSet{D_0,D_1,D_2}$ in addition to integration.

\par
Properties~\eqref{eq_twisted_zone_contained} and~\eqref{eq_twisted_zone_intersection_empty} follows directly from the $\Iverson{Y_i\subseteq{}D_i}$ terms.
Property~\eqref{eq_twisted_zone_dominate} follows from the fact that $\HSPoiThin{D_1,C_1}$ and $\HSPoiThin{D_2,C_2}$ are dominating couplings~\eqref{eq_thin_subset}.
Property~\eqref{eq_twisted_zone_connect} follows trivially from the definition of $D_1$ and $D_2$ and~\eqref{eq_twisted_zone_contained}.
\end{proof}

\begin{proposition}
\label{prop_twisted_rec}
\begin{subequations}
\label{eq_twisted_rec}
The boundary conditions of $\TwistedRec{B,C_1,C_2}$ may be restricted to $\Configs{\Ring{B}}$.
The coupling $\TwistedRec{B,C_1,C_2}$ is well-defined and jointly measurable in $(C_1,C_2)$.
Its marginals are
\begin{equation}
\label{eq_twisted_rec_hs}
 \forall\,1\le{}i\le{}2:\qquad
 \TwistedRec{B,C_1,C_2}(\xi^i=\D{Y})
 = \HardSphere{B,C_i}(\D{Y})
 \,,
\end{equation}
\begin{equation}
\label{eq_twisted_rec_poi}
 \TwistedRec{B,C_1,C_2}(\xi^3=\D{Y})
 = \Poisson{B}(\D{Y})
 \,.
\end{equation}
It has the crucial properties
\begin{equation}
\label{eq_twisted_rec_union_dominate}
 \TwistedRec{B,C_1,C_2}(\xi^1\cup{}\xi^2\subseteq\xi^3) = 1
 \,,
\end{equation}
and
\begin{equation}
\label{eq_twisted_rec_connect}
 \TwistedRec{B,C_1,C_2}((\xi^1,\xi^2)\in\DisagreementConnection{B}{C_1,C_2}) = 1
 \,.
\end{equation}
\end{subequations}
\end{proposition}

\begin{proof}
\par
The first point is to check the termination of the recursion in~\eqref{eq_twisted_def_rec}.
For $B\in\BoundedBorel$, let $\HSSize{B}:=\sup\DescSet{\Cardinality{C}}{C\in\Configs{B},\IsHS{C}=1}$.
Let $\tau:=\Diameter{B}+\Radius$.
For all $x,y\in{}B$, $\Sphere{x}$ is contained in the cube $y+[-\tau,\tau]^\Dimension$.
Putting spheres of radius $\Radius$ on a $\Radius/\sqrt{\Dimension}$ spaced $\Dimension$-dimensional integer grid within this cube covers the cube and implies that $\HSSize{B}\le(2\sqrt{\Dimension}\tau/\Radius)^\Dimension$.
Assume that there have been $n$ recursion steps~\eqref{eq_twisted_def_rec}.
This implies that there is a sequence $(x_1,\dotsc,x_n)$ of points of $B$ such that
\begin{gather*}
 \forall\,1\le{}i\le{}n-1:
 \quad{}\Distance{x_i-x_{i+1}}\le\Radius
 \,
 \\
 \forall\,1\le{}i<j\le{}n\text{ with }j-i\ge{}2:
 \quad\Distance{x_i-x_j}>\Radius
 \,.
\end{gather*}
It follows that $\IsHS{\DescSet{x_i}{1\le{}i\le{}n,i\text{ odd}}}=1$ and $\lfloor{}(n+1)/2\rfloor\le\HSSize{B}$.
Hence, the recursion terminates after at most $2\HSSize{B}$ steps.

\par
The freedom to restrict the boundary conditions to $\Configs{\Ring{B}}$ follows from the same property of $\HSPoiThin{B,.}$ in Theorem~\ref{thm_thin}, $\TwistedZone{B\setminus{}D,.,.}$ in Proposition~\ref{prop_twisted_zone} and itself.
The measurability in the boundary conditions follow from the same properties of the law $\HSPoiThin{B,.}$ in Theorem~\ref{thm_thin}, $\TwistedZone{B\setminus{}D,.,.}$ in Proposition~\ref{prop_twisted_zone}, the other measurable indicator terms in~\eqref{eq_twisted_def_zone} and itself.

\par
The marginals~\eqref{eq_twisted_rec_hs} and~\eqref{eq_twisted_rec_poi} follow directly by integrating out over the marginals and~\eqref{eq_twisted_zone_hs} and~\eqref{eq_twisted_zone_poi} respectively.
The proof of~\eqref{eq_twisted_rec_hs} use the DLR condition~\eqref{eq_dlr}.

\par
The property~\eqref{eq_twisted_rec_union_dominate} follows directly from the $\Iverson{Y_1=Y_2}$ identification and the dominating property of $\HSPoiThin{B,\emptyset}$~\eqref{eq_thin_subset} in the $D=\emptyset$ case.
In the $D\not=\emptyset$ case, it follows from~\eqref{eq_twisted_zone_contained},~\eqref{eq_twisted_zone_intersection_empty} and~\eqref{eq_twisted_zone_dominate} and itself recursively.

\par
Property~\eqref{eq_twisted_rec_connect} is trivial in the $D=\emptyset$ case and follows from~\eqref{eq_twisted_zone_connect} and itself recursively in the $D\not=\emptyset$ case.
Step $(\star)$ of the following proof of~\eqref{eq_twisted_rec_connect} in the non-trivial $D\not=\emptyset$ case demonstrates the need for the measurability of the coupling in the boundary conditions.
\begin{align*}
 &\TwistedRec{B,C_1,C_2}((\xi^1,\xi^2)\in\DisagreementConnection{B}{C_1,C_2})
 \\\stackrel{}{\ge}{}
 &\TwistedRec{B,C_1,C_2}
  ((\xi^1_{D},\xi^2_{D})\in\DisagreementConnection{D}{C_1,C_2}
  ,(\xi^1\setminus{}D,\xi^2\setminus{}D)
   \in\DisagreementConnection{B\setminus{}D}{\xi^1_{D},\xi^2_{D}}
  )
 \\\stackrel{\eqref{eq_twisted_def_rec}}{=}{}
 &\int_{\DisagreementConnection{D}{C_1,C_2}}
  \int_{\DisagreementConnection{B\setminus{}D}{Y_1,Y_2}}
  \TwistedRec{B\setminus{}D,Y_1,Y_2}((\xi^1,\xi^2)=\D{Z})
  \TwistedZone{B,C_1,C_2}((\xi^1,\xi^2)=\D{Y})
 \\\stackrel{(\star)}{=}{}
 &\int_{\DisagreementConnection{D}{C_1,C_2}}
  \TwistedRec{B\setminus{}D,Y_1,Y_2}(
   (\xi^1,\xi^2)\in\DisagreementConnection{B\setminus{}D}{Y_1,Y_2}
  )
  \TwistedZone{B,C_1,C_2}((\xi^1,\xi^2)=\D{Y})
 \\\stackrel{\eqref{eq_twisted_rec_connect}}{=}{}
 &\TwistedZone{B,C_1,C_2}((\xi^1,\xi^2)\in\DisagreementConnection{D}{C_1,C_2})
 \\\stackrel{\eqref{eq_twisted_zone_connect}}{=}{}
 &1
 \,.\qedhere
\end{align*}

\end{proof}

\section*{Acknowledgements}

I thank Marie-Colette van Lieshout, Jacob van den Berg, Ronald Meester and Erik Broman for helpful discussions surrounding this topic.
Thanks to Michael Klatt for helping me with Inkscape.
I thank the anonymous reviewers for their helpful comments.

{}
\Support{}

\bibliographystyle{plain}
\bibliography{ref}

\end{document}